\input amstex.tex

\input amsppt.sty

\TagsAsMath

\magnification=1200

\hsize=5.0in\vsize=7.0in

\hoffset=0.2in\voffset=0cm

\nonstopmode

\document

\def\N{ \Bbb N}

\def\R{ \Bbb R}

\def\Z{ \Bbb R}

\def\la{\langle}
\def\ra{\rangle}

\input amstex.tex
\input amsppt.sty
\TagsAsMath \NoRunningHeads \magnification=1200
\hsize=5.0in\vsize=7.0in \hoffset=0.2in\voffset=0cm \nonstopmode

\document

\topmatter

\title{On asymptotic stability   of  ground states of NLS with
a finite bands periodic potential in 1D }
\endtitle

\author
Scipio Cuccagna  \,and \,  Nicola Visciglia
\endauthor

\address
DISMI University of Modena and Reggio Emilia, via Amendola 2,
Padiglione Morselli, Reggio Emilia 42100 Italy\endaddress \email
cuccagna.scipio\@unimore.it \endemail

\address
Dipartimento di Matematica ``L. Tonelli'', University of Pisa, Largo
B. Pontecorvo 5, 56127 Pisa, Italy\endaddress \email
viscigli\@dm.unipi.it \endemail

\abstract We consider a nonlinear Schr\"{o}dinger equation $$iu_t
-h_0u  +  \beta ( |u|^2 )u=0 \, , \, (t,x)\in \Bbb R\times \Bbb R$$
with $h_0= -\frac {d^2}{dx^2} +P(x)$ a Schr\"{o}dinger operator with
finitely many spectral bands. We assume the existence of an
orbitally stable family of ground states. Exploiting dispersive
estimates in \cite{C2,CV} and following the argument in \cite{C1} we
prove that under appropriate hypotheses the ground states are
asymptotically stable.This paper is a slightly extended version of
the paper to be published on the Trans. AMS.

\endabstract

\endtopmatter

 \head \S 1 Introduction \endhead
Consider $h_0= -\frac {d^2}{dx^2} +P(x)$ with $P(x)\equiv P(x+1)$ a
smooth periodic potential   of period 1 in  $x\in \Bbb R$ and
consider the NLS
$$iu_t -  h_0u  +  \beta ( |u|^2 )u=0 \, ,\, (t,x)
\in \Bbb R \times \Bbb R   .\tag 1.1$$ In this paper we study  the
asymptotic stability of standing waves $u(t,x)=e^{it\omega }\phi _
{\omega }(x) $ with $\phi _\omega (x)>0$ for all $x$.
 We  consider the following hypotheses.
We  consider the following hypotheses.

{\item {(H1)}} The spectrum $\sigma (h_0)$ is formed by finitely
many bands. We choose $h_0$ so that $\inf \sigma (h_0)=0.$

{\item {(H2)}} $\beta (0)=\beta '(0)=\beta ''(0)=0$, $\beta\in
C^\infty(\R,\R)$.

{\item {(H3)}} There exists a $p\in(1,\infty)$ such that for every
$k=0,1$,
$$\left| \frac{d^k}{dv^k}\beta(v^2)\right|\lesssim
|v|^{p-k-1} \quad\text{if $|v|\ge 1$}.$$

 {\item {(H4)}} There exists
an open interval $\Cal{O}\subseteq (0,+\infty )$ such that $
    h_0u-\beta(|u|^2)u=-\omega u
$ admits a $C^1$-family of ground states $\phi _ {\omega }(x)$ for
$\omega\in\Cal{O}$.

{\item {(H5)}} $ \frac d {d\omega } \| \phi _ {\omega
}\|^2_{L^2(\Bbb R )}>0$ for $\omega\in\Cal{O}$.

{\item {(H6)}} Set $  L_+ =h_0 + \omega -\beta (\phi ^2_\omega   )
-2\phi ^2 _\omega \beta ^\prime (\phi ^2_\omega ) $. Then $\ker (L_+
)=0$ and $L_+$ has exactly one negative eigenvalue, which is simple.

\noindent By \cite{ShS} the map $ \omega \to \phi _\omega \in
H^1(\Bbb R )$   is $C^2$ and by \cite{We1,GSS} (H5) yields orbital
stability of the ground state $e^{i\omega t} \phi _ {\omega } (x)$.
Here we investigate asymptotic stability. We need some additional
hypotheses. In \S 2 we introduce the linearization $H_\omega$. By
standard arguments $ \sigma _e(H_\omega )=\cup _{\pm} \pm (\omega
+\sigma (h_0))$. In particular $ \sigma _e(H_\omega )$ is a, finite
by (H1), union of intervals. We call thresholds of $H_\omega $ the
numbers in $\Bbb R$ which are   extremes of these intervals.

{\item {(H7)}} Let $H_\omega$ be the linearized operator around
$e^{it\omega}\phi_\omega$, see (1.3). Then $H_\omega$ has a certain
number of simple positive eigenvalues with for any $j$:  $  \lambda
_j(\omega )\not \in
   \sigma _e(H_\omega )$; $2\lambda _j(\omega ) \in
\sigma _e(H_\omega )$ and is not a threshold. $H_\omega $ does not
have other eigenvalues and each threshold of $H_\omega $ is not a
resonance for $H_\omega $.

{\item {(H8)}} For   multi indexes  $m=(m_1,m_2,...)$ and
$n=(n_1,...)$, setting $\lambda (\omega )=(\lambda _1(\omega ),...)$
and $(m-n)\cdot \lambda =\sum (m_{j}-n_{j})  \lambda _j  $, we have
the following two non resonance hypotheses: {\item {(i)}}
$(m-n)\cdot \lambda (\omega )=0$ implies $m=n$ if $|m|\le 3$ and
$|n|\le 3$;  {\item {(ii)}} $(m-n)\cdot \lambda (\omega ) $ is not a
threshold of $H_\omega $ for all $(m,n)$ with $|m|+|n|\le 3$.

{\item {(H9)}} We assume the non degeneracy  Hypothesis 5.4.

\proclaim{Theorem 1.1}    Let $\omega_0\in\Cal{O}$ and
$\phi_{\omega_0}(x)$ be a ground state. Let $u(t,x)$ be a solution
of (1.1). Assume (H1)--(H9). Then, there exist an $\epsilon_0>0$ and
a $C>0$ such that if $
\inf_{\gamma\in[0,2\pi]}\|u_0-e^{i\gamma}\phi_ {\omega }
\|_{H^1(\Bbb R)}<\epsilon <\epsilon _0,$ then there exist
$\omega_+\in \Cal{O}$, $\theta\in C^1(\R;\Bbb R)$ and $  u_+\in
 H^1(\Bbb R) $    with $\| u_+\| _{H^1(\Bbb R)} \le C\epsilon  $ such that
$$\aligned &
\lim_{t\to\infty}\|u(t,\cdot)-e^{i\theta(t)}\phi_{\omega_{+}}
-e^{-ith_0}u_+\|_{H^1}=0 .
 \endaligned $$

\endproclaim

Theorem 1.1 is a transposition of more general results on the NLS
with $P(x)\equiv 0$ proved for 1D in \cite{C1} and for dimensions 2
and higher in \cite{CT,CM}. There is now a substantial literature on
asymptotic stability of ground states of the NLS, which starts with
work by Soffer and Weinstein \cite{SW1-2} and Buslaev and Perelman
\cite{BP1-2} in the early 90's, see also \cite{PW,Wd1,SW3}. In this
decade we have also further work in
\cite{C4-5,TY1-3,BS,SW4,GNT,GS1-2,M1-2}. For some results on
multisolitons see  \cite{P,RSS}. We highlight in particular
 the breakthrough work \cite{GS2},  later refined
 and extended in \cite{C1,CT,CM}. The present
 paper and the above references    focus on asymptotic stability of spatially localized standing
 waves of equation (1.1). For a discussion of
 linear and orbital stability of spatially
 periodic or antiperiodic standing waves see \cite{BR,GH}
 and references therein.

The framework for Theorem 1.1  is, by now,  classical. After
breaking  canonically solutions $u(t,x)$ into a ground state plus a
reminder, see \S 2, we consider an appropriate generalized NLS for
the reminder $R$, see (2.1).      $R$ can be decomposed using an
appropriate time varying frame associated to the spectral
decomposition of the linearization $H_{\omega (t)}$. Even though the
spectrum of $iH_{\omega (t)}$ is on the imaginary axis, the
continuous  spectrum can be thought as stable spectrum. The discrete
spectrum corresponds to central directions. As a general rule the
main difficulty consists in showing that the   discrete modes of $R$
decay to 0, that is, in some sense, showing stability in the central
manifold. The stabilization mechanism is known in various special
situations and, in our opinion, is fundamentally well understood,
although there are theorems yet to be proved. In particular, in the
references listed at the end, the Hamiltonian structure of the NLS
is not sufficiently exploited. We refer to \cite{CM,C6} for further
discussion on this theme. In this respect, the present paper
presents no novelties. With  hypotheses (H6-8) we consider a special
case where the   Hamiltonian structure of the NLS is  exploited but
which are already known in the literature. The novelty  in this
paper is that, by exploiting dispersive estimates in \cite{C2} on
$e^{ith_0}$ and extending them to $e^{itH_\omega}$, we are able to
develop the theory of radiation in the novel context involving an
NLS with a linear  periodic  potential. So far the literature on
asymptotic stability has treated only cases $h_0=-\partial
_x^2+P(x)$ with $P(x)$ either 0 or a short range potential, except
for the special case on Lam\'e potentials in \cite{C3}. We
generalize considerably \cite{C3}, in particular by considering
finite energy solutions in (1.1) without imposing further decay
conditions at infinity for the initial datum $u_0(x)$, and by easing
considerably the hypotheses on $\beta (|u|^2)$. Furthermore,
\cite{C3} treats only 2 bands potentials while in this paper $h_0$
can have any finite number of bands. To prove Theorem 1.1 we need to
use as \cite{M1-2} smoothing estimates, see also the elaborations
\cite{C1,CT}. To prove Stricharz estimates with Birman-Solomjak
  rather than Lebesgue spaces, we need to adapt a well
known lemma due to Christ and Kiselev \cite{CK}, see Lemma 3.1
\cite{SmS}.

It would be interesting to consider the case when $\sigma (h_0)$ has
infinitely many bands, but in that case $H_\omega$ can have
infinitely many eigenvalues, a situation we
  avoid by assuming that $\sigma (h_0)$ has  finitely many bands.
  Theorem 1.1 can be generalized, see \cite{CM}, assuming:

{\item {(H7')}} $H_\omega$ has a certain number of
  simple positive eigenvalues such that for any $j$ there
  an $N_j\in \Bbb N$ such that:

   $\ell \lambda _j(\omega )\not \in \sigma _e(H_\omega )$
for all $\ell \in \Bbb Z$ with $0 \le \ell \le N_j$; $(N_j+1)\lambda
_j(\omega )\in \sigma _e(H_\omega )$ but is not a threshold of
$H_\omega$. $H_\omega $ does not have other eigenvalues and each
threshold of $H_\omega$ is not a resonance. We set $N=\max _j N_j$.

{\item {(H8')}} For   multi indexes  $m=(m_1,m_2,...)$ and
$n=(n_1,...)$, setting $\lambda (\omega )=(\lambda _1(\omega ),...)$
and $(m-n)\cdot \lambda =\sum (m_{j}-n_{j})  \lambda _j  $, we have
the following two non resonance hypotheses: {\item {(i)}}
$(m-n)\cdot \lambda (\omega )=0$ implies $m=n$ if $|m|\le N+2$ and
$|n|\le N+2$;  {\item {(ii)}} $(m-n)\cdot \lambda (\omega ) $ is not
a threshold of $H_\omega$. $H_\omega $ for all $(m,n)$ with
$|m|+|n|\le N+2$.

{\item {(H9')}} There is an appropriate substitute of  Hypothesis
4.4 and we replace Lemma 4.3 with an appropriate hypothesis, see
Hypothesis 4.4 \cite{C6}.

\medskip
We remark that \cite{C1} revises \cite{C7}  but is quite different
from \cite{C7}. In \cite{C7}  Lemma 5.4 is wrong. Notice that the
material  in \S 5 \cite{C7} can be saved proceeding as in the
present paper.

\medskip

  {\bf Notation}.    For $s\in \Bbb R$,  $\| u\|_{  L^{p,s
}_x}:= \| \langle x \rangle ^su\|_{L^{p }_x} $. In particular $L^{p
}_x=L^{p,0 }_x$. We set $ \langle f ,
  g\rangle =\int _{\Bbb R}{^tf(x)} {\overline{g(x)}} dx,$  with $f(x)$ and $g(x)$ column vectors, $^tA$ the transpose and
  $\overline{g}$ the complex conjugate of $g$.
Sometimes we write $ \langle f ,
  g\rangle _x$  or $ \langle f ,
  g\rangle _t$ to emphasize the variable of integration.
  We write $ \langle f ,
  g\rangle _{sx}=\int _{\Bbb R^2}{^tf(s,x)} {\overline{g(s,x)}} dsdx.$
Given $x\in \R$ set $x^+=x\vee 0$ and $x^-=(-x)\vee 0$.
$R_H(z)=(H-z)^{-1}$.   $W^{k,p}(\Bbb R)$ is the space of tempered
distributions $f(x)$ such that $  (1-\partial _x^2) ^{k/2}f\in
L^p(\Bbb R)$. $\Bbb R^-=(-\infty ,0)$,  $\Bbb R^+=(0,  \infty  )$.
Given two  functions $f(x)$ and $g(x)$ with  values in $\Bbb C^2$,
whose elements are column vectors, then their Wronskian is the
scalar $W[f,g](x)={^tf}'(x) g(x)-{^tf} (x) g'(x).$ A pair $(r,p)$ of
numbers is said to be admissible if $(r,p) \in [4, \infty]\times [2,
\infty]$ with $2/r=1/2-1/p$. We set $\overrightarrow{e}_1 = {^t
(1,0)} $  and $\overrightarrow{e}_2 = {^t (0,1)} $. Set $\text{diag
}(a,b)$ for the diagonal $2\times2$ matrix with $(a,b)$ on the
diagonal. Given a metric space $X$,   $k\in X$ and $S\subseteq X$,
we denote the distance of $k$ from $S$ by $\text{dist}(k,S).$  For a
pair of Banach spaces $X$ and $Y$ we denote by $B(X, Y)$ the Banach
space of bounded linear operators from $X$ to $Y$. $\Cal S (\Bbb
R^n)$ is the space of Schwarz functions.

\head \S 2 Linearization, modulation and set up \endhead

We will use  the  following classical result,  \cite{We1,GSS1-2}:

\proclaim{Theorem 2.1}   Suppose that $e^{i\omega t} \phi _ {\omega
} (x)$ satisfies (H1-6) and $\phi _ {\omega } (x)>0$ for all $x$.
Then $\exists \, \epsilon >0$ and a $A_0(\omega )>0$ such that for
any  $\| u(0,x) - \phi _ {\omega } \| _{H^1_x}<\epsilon $ we have
for the corresponding solution $ \inf \{ \|  u(t,x) -e^{i \gamma
}\phi _ {\omega } (x ) \| _{H^1_x} : \gamma \in \Bbb R   \}  <
A_0(\omega ) \epsilon . $
\endproclaim      Setting $ u(t,x) = e^{i \Theta (t)} (\phi _{\omega
(t)} (x)+ r(t,x)) \, , \, \Theta (t)= \int _0^t\omega (s) ds +\gamma
(t)  $ we get

$$\aligned &
  i r_t  =
 h_0 r  +\omega (t) r-
\beta ( \phi _{\omega (t)} ^2 )r -\beta ^\prime ( \phi _{\omega (t)}
^2 )\phi _{\omega (t)} ^2 r \\&-
 \beta ^\prime ( \phi _{\omega (t)} ^2 )
\phi _{\omega (t)} ^2  \overline{  r }+ \dot \gamma (t) \phi
_{\omega (t)} -i\dot \omega (t)
\partial _\omega \phi   _{\omega (t)}
+ \dot \gamma (t) r
  +
  n(r,\overline{r} ) \endaligned
$$
for $n(r,\overline{r} )=O(r^2)$, $\overline{n(r,\overline{r}
)}=n(\overline{r}, {r}  )$. {Set }  $^tR= ( r
  ,
 \overline{ r}) $, $ ^t\Phi _\omega
=( \phi _{\omega } , \phi _{\omega }) $.   Given $ \sigma _1=\left[
\matrix  0 &
1  \\
1 & 0
 \endmatrix \right]   , $ $
\sigma _2=\left[ \matrix  0 &
i  \\
-i & 0
 \endmatrix \right]  , $ $
\sigma _3=\left[ \matrix  1 &
0  \\
0 & -1
 \endmatrix \right]  , $
 we write

$$ H_\omega =\sigma _3 \left [ h_0 + \omega   -  \beta
(\phi ^2 _{\omega  }) - \beta ^\prime (\phi ^2 _{\omega  })\phi ^2
_{\omega  } \right ] +i \sigma _2 \beta ^\prime (\phi ^2 _{\omega
})\phi ^2 _{\omega  }. $$ The equation (1.1) is Hamiltonian, because
there is a real valued function $F(|u|^2)$ with $\beta (|u|^2)
u=\partial _{\overline{u}}( F(|u|^2))$ and $F(0)=0.$ Then   $-
n(r,\overline{r} )=\partial _{\overline{r}} G(R )$ with $G(R )$ real
valued, $G(0)=0$. For $^t G'(R)=( G_{r},G_{\overline{r}})$ we
rewrite the above equation as
$$i  R _t =H _{\omega}   R +\sigma _3 \dot \gamma   R
+\sigma _3 \dot \gamma \Phi - i \dot \omega \partial _\omega \Phi
+\sigma _3 \sigma _1 G'(R).\tag 2.1
$$
Set
  $H_0(\omega )=\sigma _3(h_0 +\omega )$
and $V(\omega )=H_\omega -  H_0(\omega ).$ The essential spectrum is
$$\sigma _e =\sigma _e (H_\omega )
=\sigma _e (H_0(\omega ) ) =  \cup  _{\pm } \pm (\sigma (h_0)+\omega
) .$$ 0 is an isolated eigenvalue. Given an operator $L$ we set
$N_g(L)= \cup _{j\ge 1} N(L^j)$ and $N(L)=\ker L$. \cite{We2}
 implies that,
 if $\{ \cdot \}$ means span,
$N_g(H^\ast _\omega )=\{ \Phi , \sigma _3\partial _\omega  \Phi   \}
$.

For each $j$ we consider a generator $\xi _j\in \ker (H_\omega
-\lambda _j)$ such that $\langle \xi _j, \sigma _3 \xi _j \rangle =
1$.   We  expand $R(t)\in N_g^\perp (H^\ast _{\omega (t)}) $ into

$$    R (t) =(z\cdot \xi + \bar z\cdot \sigma _1 \xi ) + f(t)   \in
\big [ \sum _{j,\pm  }\ker (H  _{\omega (t)}\mp \lambda _j(\omega
(t)))\big ] \oplus L_c^2(H_{\omega (t)})  . \tag 2.2 $$
Correspondingly we   express  (2.1) as

$$\aligned   i\dot z _j\xi _j-\lambda _j(\omega ) z_j\xi _j&= P_{\ker (H  _\omega -\lambda _j)}
 ( \dot \gamma (\omega , R) \sigma _3R+\sigma _3 \sigma _1 G'(R)\\&
 -iz_j\dot \omega (\omega , R)
 \partial _\omega \xi _j +i\dot \omega (\omega , R) \partial _\omega P_{\ker (H  _\omega -\lambda _j)}R)
  \\  iP_c(H_\omega ) \dot f -H_\omega f &= P_c(H_\omega )
(  \dot \gamma (\omega , R) \sigma _3R+\sigma _3 \sigma _1
G'(R)+i\dot \omega (\omega , R) \partial _\omega P_{c} (H  _\omega
)R).\endaligned \tag 2.3$$

 \head \S 3 Spacetime estimates for $H_\omega $ \endhead
We list a number of linear estimates needed in the stability
argument.

 \proclaim{Lemma 3.1 (dispersive estimates)}
There exists a projection
 $\pi _{\omega} :L^2_x \to L^2_x$ with $[\pi _{\omega} ,H _{\omega}  ] =0$,
  such that if we set
  $\Cal {U}(t )= \pi _{\omega}e^{itH_{\omega}}P_c(H_{\omega})$
  and
  $\Cal {V}(t )\equiv (1-\pi _{\omega})
  e^{itH_{\omega}}P_c(H_{\omega})$,   then
we have for $C(\omega )>0$ semicontinuous in $\omega $ and for all
$t\in \Bbb R$
$$\| \Cal {U}(t ) \| _{B(L^1_x, L^\infty _x)}\le C \langle t \rangle
^{-\frac 13}\text{ and } \| \Cal {V}(t ) \| _{B(L^1_x, L^\infty
_x)}\le C | t | ^{-\frac 12}.$$

\endproclaim
The proof is sketched in \S 9.

For every $1\leq p, q\leq \infty$ we introduce the Birman-Solomjak
spaces
$$
\ell ^p({   \Bbb Z}, L^q_t[n,n+1])\equiv \left \{f\in
L^q_{loc}({\Bbb R}) \hbox{ s.t. } \{\|f\|_{L^q[n, n+1]} \}_{n\in
{\Bbb Z}} \in \ell ^p({\Bbb Z})\right \}, $$  endowed with the
  norms
$$\aligned &    \| f \|_{\ell ^p({  \Bbb Z}, L^q_t[n,n+1])} ^p\equiv
\sum_{n\in \Bbb Z}   \| f \| _{L^q_t[n,n+1] }^{p}   \quad  \forall
\quad  1\le p<\infty \text{ and } 1\leq q \leq \infty  \\ & \| f \|
_{\ell ^\infty ({ \Bbb Z}, L^q_t[n,n+1])} \equiv \sup _{n\in \Bbb Z}
\| f \| _{L^q_t[n,n+1] }.
\endaligned $$
Then we have:

\proclaim{Lemma 3.2 (Strichartz estimates)} There exists a constant
$C=C(\omega )$ upper semicontinuous in $\omega $ such that for any
$k\in [0,2]$ and for every admissible pair $(r,p) $   we have: {
$$\align &
  \| \Cal U(t) f \|    _{\ell ^{\frac 32 r}( \Bbb {Z},L^{\infty}_{t}([n,n+1],W ^{k,p}_x))} \leq
C \|f\|_{H^k_x}\tag 1\\& \| \Cal V(t) f \|
  _{ L^{r}_{t}(\Bbb R ,W ^{k,p}_x )} \leq
C \|f\|_{H^k_x}.\tag 2\endalign
$$
For any two admissible pairs $(r_1,p_1), (r_2,p_2) $   we have:
$$ \align  & \left \| \int _{0}^{t}\Cal U(t-s) g(s )ds\right \| _{\ell ^{\frac{3}{2}r_1}( \Bbb
{Z},L^{\infty}_{t}([n,n+1],W ^{k,p_1}_x)} \tag 3\\&
\le C  \| g \| _{\ell ^{\left ( \frac{3}{2}r_2\right ) '}( \Bbb
{Z},L^{1}_{t}([n,n+1],W ^{k,p_2'}_x)}
\\&
\left \| \int _{0}^{t}\Cal V(t-s) g(s )ds\right \| _{
L^{r_1}_{t}(\Bbb R,W ^{k,p_1}_x)}    \le C  \| g \| _{
L^{r_2'}_{t}(\Bbb R,W ^{k,p_2'}_x)} .\tag 4
\endalign
$$\endproclaim
Lemma 3.2 is a consequence of Lemma 3.1 and an adapted $TT^\ast$
argument, see Lemma 3.1 \cite{C1} and see \S 9 for the proof.

\medskip

In \S 7 we prove the following Kato smoothness result:

\proclaim{Lemma 3.3} Fix $\tau >3/2$, then:

{\item {(1)}} there exists $C=C(\tau ,\omega )$, upper
semicontinuous in $\omega $ such that  for any $\varepsilon \neq 0$

$$\| R_{H_\omega }(\lambda +i\varepsilon )P_c(H_\omega )
u\| _{L^2_\lambda L^{2,-\tau }_x}\le  C \| u\| _{L^2_x}; $$ {\item
{(2)}} for any $u\in L^{2, \tau }_x $ the following limits exist
$$ \lim _{\epsilon \searrow 0}R_{H_\omega }(\lambda \pm i\varepsilon )
u= R_{H_\omega }^\pm  (\lambda ) u  \text{ in $C^0(\sigma
_e(H_\omega ),L^{2, -\tau }_x)$};$${\item {(3)}} we have
$$
  \|   R_{H_\omega }^\pm  (\lambda
)P_c(H_\omega )   \| _{B( L^{2,\tau }_x, L^{2,-\tau }_x)} < C
\langle \lambda \rangle ^{-\frac{1}{2}};$$ {\item {(4)}} given any
$u\in L^{2, \tau }_x $ we have
$$P_c(H_\omega )u=\frac{1}{2\pi i}\int _{\sigma _e(H_\omega )}
(R_{H_\omega }^{+}(\lambda  )-R_{H_\omega }^{-}(\lambda  ))  u\,
d\lambda .$$
\endproclaim
Claim (1) is proved as inequality (4) in Proposition 7.1. Claim (2)
is proved as formula (12) in  Proposition 7.1. Claim (3) is a
consequence of the proof in  Proposition 7.1. The limiting
absorption principle Claim (4) is Lemma 7.6 below.

 \proclaim{Lemma 3.4} For any $k$ and  $\tau >3/2$  $\exists$
  $C=C(\tau ,k,\omega )$ upper
semicontinuous in $\omega $ such that:
   {\item {(a)}}
  for any $f\in \Cal S(\R)$,
$$\align &
\| e^{-itH_{\omega }}P_c(H_\omega )f\| _{L_{  t}^2 H_x^{k, -\tau}}
\le
 C\|f\|_{H ^{k}_x};
\endalign $$
  {\item {(b)}}
  for any $g(t,x)\in
 \Cal S (\R^2)$
$$ \left\|\int_\Bbb R e^{itH_{\omega }}
P_c(H_\omega )g(t,\cdot)dt\right\|_{H^k_x} \le C\| g\|_{L_{  t}^2
H_x^{k, \tau}}.
$$\endproclaim
  {\it Proof.} By Proposition 7.1 below, it is enough to prove
  Lemma 3.4, as well as Lemmas 3.5-6 below, for $k=0$.
   (a) implies (b) by duality:  $$\aligned & |\langle  f, \sigma _3\int_\R e^{itH_{\omega }}P_c( {\omega })g(t )dt
\rangle _{x}|= |\langle \la x\ra^{-\tau}e^{-itH_{\omega }}P_c (
H_{\omega })f, \sigma _3\la x\ra^{ \tau}g \rangle _{t,x}|\\& \le \|
e^{-itH_{\omega }}P_c( H_{\omega })f\|_{L_{  t}^2 L_x^{2, -\tau}}\|
g\|_{L_{  t}^2 L_x^{2, \tau}} \le \| f\|_{L_x^2}\| g\|_{L_{ t}^2
L_x^{2, \tau}}.
\endaligned
$$ We now prove (a) for $k=0$.
Let $g(t,x) \in \Cal S(\Bbb R^2) $ with $g(t)=P_c(H_\omega )g(t)$. Then
$$\aligned & \la  e^{-itH_\omega } f,\sigma _3
g\ra _{t,x}=  \frac{1}{\sqrt{2\pi}i}\int_\R
  \left\la
 (R_{H_\omega }^{+}(\lambda  )-R_{H_\omega
}^{-}(\lambda  ))   f,\sigma _3 \overline{\widehat{g}}(\lambda
)\right\ra_x  d\lambda  \\& = \frac{1}{\sqrt{2\pi}i}\int_{ \sigma
_e(H_\omega )}
  \left\la
 (R_{H_\omega }^{+}(\lambda  )-R_{H_\omega
}^{-}(\lambda  ))   f,\sigma _3 \overline{\widehat{g}}(\lambda
)\right\ra_x  d\lambda  .
\endaligned $$
Then from Fubini and  Plancherel and by   (1) Lemma 3.3 we have
$$\aligned & \big |\la  e^{-itH_\omega } f,\sigma _3 g\ra _{t,x}\big
|
 \le
\| (R_{H_\omega }^{+}(\lambda  )-R_{H_\omega }^{-}(\lambda ))
f\|_{L^{2, -\tau}_xL^2_\lambda (\sigma _e(H_\omega ))}     \| {g}
\|_{L^{2, \tau}_x L^2_t  }
 .
  \endaligned $$
We have
$$\align  &\| (R_{  H _\omega   }^{+}(\lambda  )-R_{  H _\omega   }^{-}(\lambda )) f\|_{\ell ^{2,
-\tau} L^2_\lambda ({\sigma _e(  H_\omega )})}=\\& \lim
_{\varepsilon \searrow 0} \| (R_{  H _\omega  } (\lambda
+i\varepsilon )-R_{  H _\omega   } (\lambda -i\varepsilon ))
f\|_{\ell ^{2, -\tau} L^2_\lambda ({\sigma _e(  H_\omega
)})}.\endalign
$$
It is enough to prove Lemma 3.4 for $f=P_c(  H _\omega )f$. Then
  $$\| R_{  H _\omega  } (\lambda +i\varepsilon )
f\|_{\ell ^{2, -\tau} L^2_\lambda ({\sigma _e(  H_\omega )})} \le C
\| f\| _{\ell ^2}
$$
for any $\epsilon \neq 0$ and fixed $C$ follows from inequality (4)
in the proof of Proposition 8.1 below.

\proclaim{Lemma 3.5} For any $k\in \Bbb N$ and  $\tau >3/2$
$\exists$
  $C=C(\tau ,k,\omega )$
such that $\forall$ $g(t,x)\in \Cal S(\R^2)$
$$\align &  \left\|  \int_0^t e^{-i(t-s)H_{\omega
}}P_c(H_\omega )g(s,\cdot)ds\right\|_{L_{  t}^2 H_x^{k, -\tau}} \le
C\|  g\|_{L_{  t}^2 H_x^{k, \tau}}.\endalign
$$
\endproclaim
{\it Proof.}  By Plancherel and H\"older inequalities and by (3)
Lemma
 3.3   we have
$$ \aligned &
\| \int _{0}^t e^{-i(t-s)H_{\omega }}P_c(H_\omega )g(s,\cdot)ds\|_{
L_{t }^2L_{ x}^{2,-\tau }}  \\& \le  \| R_{H_\omega }^+(\lambda
)P_c (H_\omega )
  \widehat{ \chi }_{[0,+\infty )}\ast _\lambda
   \widehat{ g}(\lambda,x)\|_{L_{\lambda }^2L_{ x}^{2,-\tau } }  \\& \le
   \left\| \,
\|   R_{H_\omega }^+ (\lambda )P_c  (H_\omega ) \| _{B(
L^{2,\tau}_x, L^{2,- \tau}_x)} \|
   \widehat{ \chi }_{[0,+\infty )}
   \ast _{\lambda } \widehat{g} (\lambda,x) \|_{L_{ x}^{2, \tau }}\, \right\|_{L^2_\lambda}
\\ \le &
  \|  R_{H_\omega }^+ (\lambda
)P_c(H_\omega )   \| _{L^\infty _\lambda (\Bbb R ,B( L^{2,\tau}_x,
L^{2,-\tau}_x))}\| g\|_{L_{t }^2L_{ x}^{2, \tau } }  \le C \|
g\|_{L_{t }^2L_{ x}^{2, \tau } } .
\endaligned $$

\bigskip

\proclaim{Lemma 3.6} For any $k\in \Bbb N$ and  $\tau >3/2$
$\exists$
  $C=C(\tau ,k,\omega )$ such that $\forall$ $g(t,x)\in \Cal S(\R^2)$
$$\align  &
\left\|\int_0^t \Cal U(t-s)g(s,\cdot)ds \right\|_{  L_t^\infty
L_x^2\cap \ell ^{6}( \Bbb {Z},L^{\infty }_{t}([n,n+1],W ^{k,\infty
}_x))} \le C\|g\|_{L_t^2H_x^{k,\tau}} \tag 1\\& \left\|\int_0^t \Cal
V(t-s)g(s,\cdot)ds \right\|_{  L_t^\infty L_x^2\cap  L^{4 }_{t} W ^{k,\infty }_x) } \le C\|g\|_{L_t^2H_x^{k,\tau}}.\tag 2
\endalign$$
\endproclaim
 {\it Proof.} For $g(t,x)\in \Cal S(\Bbb R^2)$ set
$$T_1g(t)=\int _0^{+\infty}  \Cal U(t-s)g(s) ds \text{ and }
T_2g(t)=\int _0^{+\infty}  \Cal V(t-s)g(s) ds.$$ Lemma 3.4 (b)
implies $f:=\int _0^{+\infty}  e^{isH_\omega }P_c(H_\omega
)g(s)ds\in L^2_x$. So in particular $\pi _\omega f$ and $(1- \pi
_\omega ) f$ belong to $L^2_x$. For $(r,p) $ admissible we have
$$ \aligned &  \|T_1g(t)\|_{\ell ^{\frac 32 r}
( \Bbb {Z},L^{\infty}_{t}([n,n+1],L^{ p}_x))}\lesssim  \| \pi_\omega
f\|_{L^2_x} \lesssim  \|g\|_{L^2_tL^{2, \tau}_x }\\& \|T_2g(t)\|_{
L^{r}_{t} L^{ p}_x }\lesssim  \| (1-\pi_\omega ) f\|_{L^2_x}
\lesssim  \|g\|_{L^2_tL^{2, \tau}_x }.\endaligned
$$ We get
 as a direct consequence of \cite{CK}, Lemma 3.1 \cite{SmS}:
$$
\left\|\int_{0}^{t}   \Cal V (t-s) g(s)ds \right\| _{   L^{r}_{t}
L^{ p}_x } \lesssim \|g\|_{  L^2_tL^{2, \tau}_x  }.
$$
Our next claim is that
$$
\left\|\int_{0}^{t}   \Cal U (t-s) g(s)ds \right\| _{\ell ^{\frac 32
r}( \Bbb {Z},L^{\infty }_{t}([n,n+1],L^{ p}_x))} \lesssim \|g\|_{
L^2_tL^{2, \tau}_x  }.
$$
For $(r,p)=( \infty ,2)$ this follows from \cite{CK}. The case
 $(r,p)=(
4,\infty )$ follows from an extension of to Birman-Solomjak  spaces
of the result in \cite{CK} which we prove in     section 10.

\head \S 4  Nonlinear argument for system (2.3) \endhead

 We   use the multi index notation
$z^m=\prod _j z^{m_j}_j$. We  consider the Taylor expansion
$$\aligned & \sigma _3 \sigma _1 G'(R)=\sum _{   |m+n|=2}^{ 3} R_{m,n}(\omega ) z^m  \bar z^n+
\sum _{ |m + n|=1}  z^m  \bar z^n A_{m,n}(\omega ) f+ O(f^2)+\cdots
\endaligned    $$ with $R_{m,n}(\omega  ,x) $ and
$A_{m,n}(\omega ,x ) $ real vectors  and matrices exponentially
decreasing in $x$. We have
$$\aligned & A_{m,n}(\omega )= \frac{\sigma _3 \sigma _1}{m! n!}\partial _z^m
 \partial _{\overline{z}} ^n\partial _f G'(0)
  \, , \quad  R_{m,n}(\omega )= \frac{\sigma _3
\sigma _1}{m! n!}\partial _z^m  \partial _{\overline{z}} ^nG'(0)
 .\endaligned$$
  We set
$\delta _j=(\delta _{j1},\delta _{j2},...) $
 with $\delta _{jk}$ the Kronecker delta. We
have
$$ \aligned &  A_{\delta _\ell ,0}(\omega ) =
\sigma _3 \sigma _1 \partial _{z_\ell}
 \partial _f G'(0)  = \sigma _3 \sigma _1\partial
_{z_\ell }
 \partial _f G'(0)  =  \sigma _3 \sigma _1 G ^{(3)}(0)
    ( \cdot  ,\xi _\ell , P_c(H  _\omega  ))  \endaligned $$
where $G ^{(3)}(0)$ is written as a symmetric trilinear form and
where one of the vectors of the triple is $\xi _\ell$. We have

$$ \aligned & P_c(H  _\omega  )R_{ \delta _j +\delta _\ell ,0}(\omega ) =
\frac{\sigma _3 \sigma _1}{(\delta _j +\delta _\ell )!}P_c(H _\omega
^\ast )\partial _{z_j}\partial _{z_\ell }
   G'(0)  \\& = \frac{\sigma _3 \sigma _1}{(\delta _j +\delta _\ell )!}\partial
_{z_j}\partial _{z_\ell }
   G'(0) \circ P_c(H  _\omega    )=  \frac{\sigma _3 \sigma _1}{(\delta _j +\delta _\ell )!} G ^{(3)}(0)
 (\xi _j,\xi _\ell  ,  P_c(H  _\omega  ))  .\endaligned $$

For  $O_{loc}(z^n)=\sum _\ell O_{loc}(|z_\ell ^{n}|)$, (2.1) can be
expressed as

$$\aligned &if_t=\left ( H _{ \omega (t)}+P_c(H_\omega )\sigma _3 \dot \gamma \right )f
  +
\sum _{  |m+n|=2} z^m \bar z^nP_c(H_\omega )R_{m,n}(\omega )
 \\& + \sum _{  |m + n|=1} z^m \bar z^n P_c(H_\omega )A_{m,n}(\omega ) f+
O(f^2)+O_{loc}(z^3),
\endaligned \tag  4.1a
$$
 and
$$\aligned & i\dot z _j\xi _j-\lambda _j(\omega ) z_j\xi _j = P_{\ker (H  _\omega -\lambda _j)}
 ( \sum _{  |m+n|=2} z^m \bar
z^nR_{m,n}(\omega )\\& +\sum _{  |m+n|=3} z^m \bar
z^nR_{m,n}^{(1)}(\omega )
 +\sum _{  |m + n|=1} z^m \bar z^n A_{m,n}(\omega ) f)+O(f^2)+O_{loc}(z^4)
\endaligned \tag  4.1b
$$
with the $A_{m,n}(\omega )$ and $R_{m,n}(\omega )$   the same of the
expansion of $\sigma _3 \sigma _1 G'(R)$   and  with the
$R_{m,n}^{(1)}(\omega )$   real
  and exponentially decreasing vectors.
We have:

  \proclaim{Lemma 4.1}  For any $C_1>0$
$\exists$ $\epsilon (C_1)>0$ and   $C(C_1)$ such
  that if, for $0<\epsilon <\epsilon (C_1)$, we have
  $\| z_j\| _{L^\infty (0,T)}+\| z_j\| _{L^4(0,T)}^2\le C_1\epsilon $ for all $j$, then

$$\|
  f  \| _{L^2_t((0,T),H^{1,-2}_x)}+\| f\|
_{L^6_t((0,T),L^\infty _x)\cap L^\infty _t((0,T),H^1 _x)}<
C(C_1)\epsilon . $$
\endproclaim
{\it Proof.} Notice that $\sum _j\| z_j\| _{L^\infty_t} +\|
f\| _{  L^\infty _t H^1 _x}<c_0 \epsilon $ for a fixed $c_0 $
is a consequence of orbital stability, Theorem 2.1 above. In this
proof we set $P_c(\omega )=P_c(H _{\omega} )$. We split $P_c(\omega
)= P_+(\omega )+P_-(\omega )$, with  $P_\pm (\omega )$ the spectral
projections in $\Bbb R_\pm \cap \sigma _e(H_\omega )$, see Lemma
8.13 below. We now rewrite the equation
  for
  $f$.
By orbital stability, Theorem 2.1, we can fix $\omega _0$ such that
$|\omega (t)-\omega _0|=O(\epsilon )$ for all $t$. Following
\cite{BP2} we set
$$\aligned &if_t=\left \{ H_ {\omega _0}   + (\dot \gamma +\omega
-\omega _0) (P_+(\omega _0)-P_-(\omega _0))\right \} P_c(\omega _0)f
\\& +O_{loc}(\epsilon f)+
O(f^2)+O_{loc}(z^2)\\& \text{where } O_{loc}(\epsilon f)=
 (\dot \gamma +\omega -\omega _0) \left (P_c(\omega _0)\sigma _3-  (P_+(\omega _0)-P_-(\omega _0))  \right ) f
 \\&  +
   \left (  V  (\omega ) -
 V  (\omega _0) \right )   f
  +(\dot \gamma +\omega -\omega _0) \left ( P_c(\omega  )-  P_c(\omega _0)\right )\sigma _3f
  .
\endaligned
$$
with $V  (\omega )$ the localized matrix potential defined under
(2.1) and where  the notation $O_{loc}(\epsilon f)$ is justified by
the fact that  $\omega -\omega _0=O(\epsilon )$ and for any pair
$s_1,s_2\in \Bbb R$ there is $ c_{s_1,s_2} (\omega )$ upper
semicontinuous in $\omega$, see Lemma 8.13, such that for  $j=0,1$

$$   \|   P_c (\omega
) \sigma _3- (P_+ (\omega )-P_- (\omega ) ) \|  _{H ^{j,s_1}_x\to  H
^{j,s_2}_x}    \le c_{s_1,s_2} (\omega )<\infty  $$  where the case
$j=1$ follows from $j=0$ applying $H_\omega$ and interpolating.  We
have $\| f \| _{L^6_t((0,T),L^\infty _x)}\approx \| P_c(\omega _0)f
\| _{L^6_t((0,T),L^\infty _x)}$ by

$$\aligned & \| f \| _{L^6_t((0,T),L^\infty _x)}= \| P_c(\omega _0)f
\| _{L^6_t((0,T),L^\infty _x)} +\| \left ( P_c(\omega  )- P_c(\omega
_0)\right )f \| _{L^6_t((0,T),L^\infty _x)} \\& = \| P_c(\omega _0)f
\| _{L^6_t((0,T),L^\infty _x)} + O(\epsilon)\| f \|
_{L^6_t((0,T),L^\infty _x)}\endaligned $$ For $$ U_\pm (t,t')=
e^{-i(t-t')H _{\omega_0}} e^{\pm i\int _{t'}^t d\tau (\dot \gamma
(\tau )+\omega (\tau )-\omega _0) }P _{\pm}(\omega _0),\tag 4.2$$ we
have
$$\aligned & P _{\pm}(\omega _0)f(t)=  U_\pm (t,0)f(0)+\int _0^tU_\pm (t,t')
(O_{loc}(\epsilon f)+ O(f^2)+O_{loc}(z^2)) dt'.
\endaligned $$
We   set $O_{loc}(\epsilon f)+ O(f^2)+O_{loc}(z^2)=\Cal X + O(f^7).$
Let $\pi =\pi _{\omega _0}$ be the projection in Lemmas 3.1-2   and
$\pi _1=1-\pi$. By (1)--(2) Lemma 3.2 and interpolation we have for
a fixed $c_0(\omega _0)$
$$\aligned &\|   U_\pm (t,0)f(0)\| _{L^6_t((0,T),L^\infty _x)} \le
 \|   U_\pm (t,0)f(0)\| _{L^6_t L^\infty _x}\lesssim \\&
  \|   U_\pm (t,0)\pi  f(0)\| _{\ell ^6( \Bbb Z,
  L^\infty _t([n,n+1] ,L^\infty _x)) } +
  \|   U_\pm (t,0)\pi _1f(0)\| ^{\frac{2}{3}}_{L^4_t L^\infty_x}
  \\& \times  \|   U_\pm (t,0)\pi _1f(0)\| ^{\frac{1}{3}}_{L^\infty _t H^1_x}    \le c_0(\omega_0 )
  \| f(0) \| _{H^1_x} \le c_0(\omega _0)\epsilon .\endaligned $$
 Similarly, by (3)--(4) Lemma 3.2 and interpolation we have   $$\aligned & \left \| \int _0^tU_\pm (t,t') O(f^7) dt'
\right \| _{L^6_t((0,T),L^\infty _x)}\lesssim \left \| \int
_0^tU_\pm (t,t')\pi _{1} O(f^7) dt' \right \| ^{\frac{2}{3}}
_{L^4_t((0,T),L^{\infty }_x)}
\\& \times
\left \| \int _0^tU_\pm (t,t')\pi _{1} O(f^7) dt' \right \|
^{\frac{1}{3}} _{L^\infty _t((0,T),H^{1 }_x)}
\\& + \left \| \chi _{[0,T]}\int _0^tU_\pm (t,t')\pi   O(f^7)
dt' \right \| _{\ell ^6( \Bbb Z,
  L^\infty _t([n,n+1] ,L^\infty _x)) }.\endaligned $$
By Lemma 3.2 both terms in the right are bounded by $\| O(f^7)  \|
_{L^1_t((0,T),H^{1 }_x)}$. So by Lemma 3.2 and H\"older inequality
we have
$$\aligned & \left \| \int _0^tU_\pm (t,t') O(f^7) dt'
\right \| _{L^6_t((0,T),L^\infty _x)} \le c_0 \| O(f^7)  \|
_{L^1_t((0,T),H^{1 }_x)}
\\& \lesssim
\| f\| ^6 _{L^6_t((0,T),L^\infty _x) }   \| f\| _{L^\infty
_t(0,T)H^1 _x}\lesssim \epsilon ^{7}.\endaligned $$ In a similar
fashion, by Lemma 3.6  and by (1) below we have
$$ \left \| \int _0^tU_\pm (t,t')
\Cal X dt' \right \| _{L^6_t((0,T),L^\infty _x)} \le c_0(\omega ) \|
 \Cal X \| _{L^2_t((0,T),H^{1,2}_x)} \lesssim
\epsilon ^{2}.$$ We claim we have
$$ \|   \Cal X \| _{L^2_tH^{1,2}_x}\le C (\omega )
\left ( \epsilon ^2+ \|   f  \| _{L^2_tH^{1,-2}_x} ^2\right  ).\tag
1$$ Indeed $\Cal X=O_{loc}(\epsilon f)+O_{loc}(z f)+
O_{loc}(f^2)+O_{loc}(z ^2)$ with
  $$ \|  O_{loc}(\epsilon f) \| _{  H_x^{1,2} L_t^2 }
\lesssim \epsilon \| f \| _{L_t^2 H_x^{1,-2}} \, ; \quad \| O_{loc}(
z f ) \| _{  H_x^{1,2} L_t^2 } \lesssim \| z\| _{\infty} \|   f  \|
_{L_t^2 H_x^{1,-2}} ;$$
  $$ \|  O_{loc}(    f  ^2  ) \| _{ L_t^2 H_x^{1,2} }
  \lesssim   \|   f \| _{L_t^2 H_x^{1,-2} }^2 \, ; \quad \|
O_{loc}( |z |^{2} ) \| _{  H_x^{1,2} L_t^2 } \lesssim
   \epsilon  \|z  ^{2}\| _{L^2_t} \lesssim \epsilon ^2 .$$
Hence (1) is correct. Let now $f = g+h$ with
$$\aligned & i g _t=
\left \{ H _{\omega _0}   + \ell (t) (P_+(\omega _0)-P_-(\omega
_0))\right \}g+\Cal X \, , \quad  g(0)=f (0)\\&   i h _t= \left \{ H
_{\omega _0}   + \ell (t) (P_+(\omega _0)-P_-(\omega _0))\right
\}h+O(f ^7) \, , \quad  h(0)=0.\endaligned $$   Then, by Lemmas 3.4
and 3.5   we get  for a fixed $c_0$ $$ \|   g\| _{L^2_t((0,T),
H_x^{1,-2})}\lesssim c_0(C_1+1)\epsilon +O(\epsilon ^2)\le
C(C_1)\epsilon .
$$   Finally by  Lemma 3.4
$$\aligned &\int _0^T  \|
 e^{-i(t-s) H_{\omega _0}}  e^{\pm i \int _s^t\ell (\tau ) d\tau }
O(f ^7)  (s)\| _{L^2_t((0,T), H_x^{1,-2})}    \\& \lesssim \int _0^T
\| O(f ^7) (s) \| _{ H^{1}_x}ds\lesssim \epsilon  ^7 .
\endaligned $$
  This  yields
$\|   h\| _{L^2_t((0,T), H_x^{1,-2})} \lesssim \epsilon^7$ and
completes the proof of Lemma 4.1.

\bigskip
Having obtained Lemma 4.1, we rewrite  (4.1a) in the form

$$\aligned &if_t=\left \{ H_ {\omega _0}   + (\dot \gamma +\omega
-\omega _0) (P_+(\omega _0)-P_-(\omega _0))\right \}  f
\\& +
\sum _{  |m+n|=2} z^m \bar z^nP_c( \omega _{0})R_{m,n}(\omega )
  + \sum _{  |m + n|=1} z^m \bar z^n P_c( \omega _{0})A_{m,n}(\omega ) f+
\\& +
 (\dot \gamma +\omega -\omega _0) \left (P_c(\omega _0)\sigma _3-  (P_+(\omega _0)-P_-(\omega _0))  \right ) f  +
   \left (  V  (\omega ) -
 V  (\omega _0) \right )   f
\\& +(\dot \gamma +\omega -\omega _0) \left ( P_c(\omega  )-  P_c(\omega _0)\right )\sigma
_3f + O(f^2)+O_{loc}(z^3)\\& + ( P_c( \omega )-P_c( \omega _0))
\left (\sum _{ |m+n|=2} z^m \bar z^nR_{m,n}(\omega )
  + \sum _{  |m + n|=1} z^m \bar z^n  A_{m,n}(\omega ) f\right ) .
\endaligned
$$
 We then set

$$ f_{2}=f+ \sum _{  |m+n|=2} R_{H_{\omega _0}}^{+}((m-n)\cdot \lambda (\omega _0)
) P_c( H _{\omega _0})R_{m,n}    (\omega  )z^m  \bar z^n.\tag 4.3
$$
We will need below:

\proclaim{Lemma 4.2} Assume the hypotheses of Lemma 4.1.  Then for
$s>1$ sufficiently large we can decompose $f _{2}= h_1+h_2+h_3+h_4$
with: {\item {(1)}} for a fixed $c_0(\omega _0), $ $\| h _1\|
_{L^2_t L^{2,-s}_x} \le c_0(\omega _0) \| f(0)\| _{H^1}\le
c_0(\omega _0) \epsilon  ;$ {\item {(2)}} for a fixed $c_1(\omega
_0), $ $\| h _2\| _{L^2_t L^{2,-s}_x} \le c_1(\omega _0) |z(0)|^2\le
c_1(\omega _0)\epsilon ^2;$ {\item {(3)}} $\| h _3\| _{L^2_t
L^{2,-s}_x} = O(\epsilon ^2 ) ;$ {\item {(4)}}     $ \| h_4 \|
_{L^2_t L^{2,-s}_x}  =O(\epsilon ^2 ).$

\endproclaim
{\it Proof.} The proof is basically that in \S 4 \cite{CM}. We have
schematically

$$\aligned &i\partial _tP_c(H _{\omega _0})f_{2}=\left ( H _{\omega _0}   + (\dot \gamma +\omega -\omega _0)
(P_+(\omega _0) -P_-(\omega _0) )\right ) P_c(H _{\omega _0}) f_{2}
+
\\& + \sum _{|m +n|=2 } O(|z|^{3}) \
 R_{H _{\omega _0}}^+((m-n)\cdot\lambda (\omega _{0} )    R_{m,n}
(\omega _{0} )
\\&  +
  P_c(H _{\omega _0})\left (o(1) O_{ loc }(|z|^{2}) +o(1) O_{ loc }(f)+O(f^2)\right ) .  \endaligned  $$
For $h_1(0)=f(0)$ let
$$ i\partial _t (h_1+h_2)=\left ( H_{\omega _0}   + (\dot \gamma +\omega -\omega _0) (P_+-P_-)\right )
(h_1+h_2) , \quad  h_1(0)+h_2 (0)=f _{2}(0)  .$$ Then (1) follows by
Lemma 3.3 applied to $P_\pm (\omega _0)h_1(t)=U_\pm (t,0) f(0)$,
with $U_\pm (t,s)$ defined in Lemma 4.1. To get (2) we observe
  that for a constant $C=C(\Lambda , \omega  _0)$
upper semicontinuous in $\omega _0 $ and in $\Lambda
>\omega _0$ we have

$$
\|   U_\pm (t,t')R _{H_{\omega  }}^{+}(\Lambda     ) P_c g\|
_{L_x^{2,-s}} < C\langle t-t' \rangle ^{-\frac 32} \|  g \|
_{L_x^{2,s}} \, ,\, s> s _0 . \tag 5
$$
The weighted estimate (5) can be proved proceeding along the lines
of the  proof of Theorem 2.4 in pp. 135--6 \cite{C3} using the
estimates for the Jost functions in Lemma 9.1 and the representation
in Lemma 9.11. Then (5) implies (2)  by  $f _{2}(0) =\sum _{
|m+n|=2} R_{H_{\omega _0}}^{+}((m-n)\lambda (\omega _0 )) R_{m,n}
(\omega _0)z^m(0)  \bar z^n (0)  $. Next we define $h_3(0)=0$ and
$$\aligned &i\partial _t P_c(H _{\omega _0}) h_3=\left ( H_{\omega _0}   + (\dot \gamma +\omega -\omega _0)
(P_+(\omega _0) -P_-(\omega _0) )\right ) P_c(H _{\omega _0}) h_3 +
\\&
+ P_c(H _{\omega _0})\left (O(\epsilon ) O_{ loc }(|z|^{2})
+O(\epsilon ) O_{ loc }(f)+O(f^2)\right )   . \endaligned  $$  Then
(3) follows from the argument in Lemma 4.1. Finally we set
$h_4(0)=0$ and

$$\aligned &i\partial _t h_{4}=\left ( H _{\omega _0}
  + (\dot \gamma +\omega -\omega _0) (P_+(\omega _0) -P_-(\omega _0))\right ) h_{4}  \\& +
\sum _{|m +n|=2 } O(|z|^{3})
 R_{H _{\omega _0}}^+((m-n)\cdot\lambda (\omega _{0} ))    R_{m,n}
(\omega _{0} ) .\endaligned $$ Then we have $h_{4}= \sum _\pm h_{4
\pm }$ with

$$h_{4  \pm } (t)=\sum _{|m+n|=2 }\int _0^tU_\pm (t,t')
 O(|z(s)|^{3}) R_{H _{\omega _0}}^+((m-n)
 \cdot\lambda (\omega _{0} )  ) R_{m,n}
(\omega _{0} ) dt' .$$   By (5) we get $ \|h_{4 \pm } (t)  \|  _{L_x
^{2,-s}}\le C \epsilon \int _0^t \langle t-t'\rangle ^{-\frac 32}
|z(t')|^{2}  dt'$  and so  $ \| h _4\| _{L^2_{t }L_x ^{2,-s}} \le
\epsilon \| z \| _{L^{4} _{t }} ^{2}= O(\epsilon ^2   ). $ This
concludes the proof of Lemma 4.2.

\bigskip

By substitution of (4.3) in    (4.1b) we get
$$\aligned &i\dot z _j\xi _j-\lambda _j(\omega ) z_j\xi _j = P_{\ker (H  _\omega -\lambda _j)}
 ( \sum _{  |m+n|=2} z^m \bar
z^nR_{m,n}(\omega ) +\sum _{  |m+n|=3} z^m \bar
z^nR_{m,n}^{(1)}(\omega )\\& -\sum _{  |m' + n'|=1}\sum _{  |m +
n|=2} z^{m +m'}\bar z^{n+n'} A_{m',n'}(\omega ) R_{H_{\omega _0
}}^{+}((m-n)\cdot \lambda (\omega _0) ) P_c( H _{\omega })R_{m,n}
(\omega )\\&
 +\sum _{  |m + n|=1} z^m \bar z^n  A_{m,n}(\omega ) f_2+O(f^2)+O_{loc}(z^4).
\endaligned
$$
Here recall $P_{\ker (H  _\omega -\lambda _j)}=  \xi _j \langle
\quad , \sigma _3\xi _j \rangle $. By standard normal forms
arguments there exists a change of variables $ {\zeta }_j=z_j +\sum
_{|m+n|=2}^3p_{j,m,n}(\omega )z^m\overline{z}^n$ with $\Im
p_{j,m,n}=0$    for $|m+n|=2$ such that below we have $\Im a_{j,\ell
} (\omega )=0$ and

$$\aligned &i {\dot {\zeta}} _j\xi _j-\lambda _j(\omega ) \zeta_j\xi _j = \xi _j
\sum _{\ell}a_{j,\ell } (\omega )|\zeta _\ell |^2\zeta
_j+O_{loc}(\zeta f_2) +O(f^2)+O_{loc}(\zeta ^4)\\& - \sum _\ell
|\zeta _\ell |^2\zeta _jP_{\ker (H  _\omega -\lambda _j)}
A_{0,\delta _{\ell }}(\omega )R_{H_{\omega _0}}^{+}( \lambda _{\ell
}(\omega _0 ) +\lambda _j (\omega _0) ) P_c( H _{\omega })R_{\delta
_{\ell }+\delta _j,0} (\omega  )
\endaligned \tag  4.4
$$
where   $O_{loc}(\zeta ^n f_2)=\sum
 _\ell O_{loc}( \zeta _\ell ^{n}f_2)$.
Applying  $\langle \quad , \sigma _3\xi _j\rangle $ to (4.4)  we get

$$\aligned &  i {\dot {\zeta}} _j -\lambda _j(\omega ) \zeta_j  =
\sum _{\ell}a_{j,\ell } (\omega )|\zeta _\ell |^2\zeta _j+\langle
O_{loc}(\zeta f_2) +O(f^2)+O_{loc}(\zeta ^4) ,    \sigma _3\xi
_j\rangle  -\\& -   \sum _\ell |\zeta _\ell |^2\zeta _j \langle
A_{0,\delta _{\ell }}(\omega )R_{H_{\omega _0}}^{+}( \lambda _{\ell
}(\omega _0 ) +\lambda _j (\omega _0) ) P_c( H _{\omega })R_{\delta
_{\ell }+\delta _j,0} (\omega  ) , \sigma _3\xi _j\rangle .
\endaligned \tag 4.5
$$
Recall that $\Im a_{j,\ell } (\omega )=0$. Then multiplying (4.5) by
$\overline{\zeta}_j$ we get
$$\aligned &   \frac{1}{2} \frac{d}{dt}{| \zeta _j|^2}   =
-   \sum _\ell  \left ( \Gamma _{j,\ell } (\omega _0) +o(1)\right )
|\zeta _\ell |^2|\zeta _j|^2+\\& +\Im [ \overline{\zeta _j}\langle
O_{loc}(\zeta f_2) +O(f^2)+O_{loc}(\zeta ^4) ,    \sigma _3\xi
_j\rangle ]
\endaligned \tag 4.6
$$
where we use $\omega -\omega _0=O(\epsilon )$,
$$\aligned & \Gamma _{j,\ell } (\omega  )=
 \Im \langle A_{0,\delta _{\ell
}}(\omega  )R_{H_{\omega   }}^{+}( \lambda _{\ell }(\omega   )
+\lambda _j (\omega  ) ) P_c( H _{\omega  })R_{\delta _{\ell
}+\delta _j,0} (\omega   ) , \sigma _3\xi _j(\omega   )\rangle
\endaligned $$
and the continuous dependence  in $\omega $   of $A_{0,\delta _{\ell
}}(\omega )$, $R_{\delta _{\ell }+\delta _j,0} (\omega  ) $ and
$P_c( H _{\omega })$. We have:

\proclaim{Lemma 4.3}  We have for $h_{j,\ell}(\omega )=\sigma
_3\sigma _1G ^{(3)}(0)
 (\xi _j,\xi _\ell  ,  P_c(H  _\omega  ))$
$$\aligned &  \Gamma _{j,\ell } (\omega  ) =    \frac{\pi }{(\delta _j +\delta _\ell
)!}  \langle
 \delta ( H_\omega   -\lambda _j - \lambda _\ell )
   h_{j,\ell}(\omega )
  ,  \sigma _3 h_{j,\ell}(\omega ) \rangle   \ge 0.
\endaligned $$ \endproclaim
The proof is in Lemma 3.4 \cite{C6}. We assume the following, which
by $2\lambda _j(\omega  )  >\omega $ for any $j$  looks like a
generic condition:

\proclaim{Hypothesis 4.4} We suppose that  $\Gamma _{j,j } (\omega
)>0$    for any $j$ .
\endproclaim

Integrating    (4.6) in an interval $[0,t]$ we get
$$  \aligned & C_0\sum _{j} \int _0^t|\zeta _j|^4  + \sum
_j {|{\zeta _j}(t)|^2}/{2}\le \\& \le \sum _j {|{\zeta_j}(0)|^2}/{2}
+c_1  \left ( \sum _{j} \int _0^t|\zeta _j|^4\right ) ^{\frac{1}{2}}
\| f _{2}\| _{L_t^2 L^{2,-s}_x}+O(\epsilon ^3).\endaligned
$$

By  Lemma 4.2,  $\| f _{2}\| _{L_t^2 L^{2,-s}_x}\le c_2\epsilon$ for
$c_2\approx 1$. Then for a fixed $c\approx 1$
$$ \sum _{j}\int _0^t|z _j|^4  + \sum
_j {|  {z}_j (t)|^2} \le c \, \epsilon ^2.$$

\bigskip

The proof that, for $^tf (t)=(h(t),\overline{h}(t))$, $h(t)$ is
asymptotically free for $t\to \infty$, is similar to the analogous
one in \cite{CM}. For $  U_\pm (t,t') $ defined in (4.2) we have for
 $t_1<t_2$
$$\aligned & \| U(0,t_2)  f(t_2)-
U(0,t_1)  f(t_1) \| _{H^1_x}\le  \\& \le \| U_\pm (0,t')
(O_{loc}(\epsilon f)+ O(f^2)+O_{loc}(z^2))\| _{L^1((t_1,t_2),H^1_x)}
\lesssim \\&   \| f\| _{ L^2((t_1,t_2), H_x^{1,-2}} +\| f\| _{
L^6((t_1,t_2), L_x^{\infty }} +|z| _{L^4 (t_1,t_2)}\to 0
  \text{ for $t_1\to \infty$}.
\endaligned $$
Then   consider $w_+=P _{+}(\omega _0)w+P _{-}(\omega _0)w$  with
 $P _{\pm}(\omega _0)w=\lim _{t\to \infty}U_\pm (0,t )
   f(t ) .$ We have

$$ \lim _{t \to \infty } \| P_c(\omega _0)f (t) \|_{ L_x^{2,-2}}
=\lim _{t \to \infty } \| ( U_+(t,0)+U_-(t,0))P_c(\omega _0)f (0)
\|_{ L_x^{2,-2}}=0
$$
with the second equality true for any $f(0) \in H^1 $. For
$P_d(\omega )=1-P_c(\omega )$ we have

$$\aligned
\| P_c(\omega _0) {f} (t)-f (t)\|_{H^1}=& \|(P_d(\omega
)-P_d(\omega_0 ))f \|_{L^2_x}
\\ \lesssim & |\omega-\omega_0|
 \|P_c(\omega _0) {f} (t) \|_{ L_x^{2,-2}}
\lesssim \|P_c(\omega _0) {f} (t) \|_{ L_x^{2,-2}}\to 0,
\endaligned $$
as $t\to\infty$. Combining the above we have for $\theta (t)=\int
_{0}^t d\tau (\dot \gamma (\tau )+\omega (\tau ) )$

$$\lim_{t\to  \infty}\|f (t)-e^{i\left [ (t\omega
_0-\theta (t)+\theta (0) \right ] (P_+(\omega _0)- P_-(\omega _0))}
e^{-itH _{\omega _0}} w _+\|_{H^1_x}=0 .$$ Consider the strong limit $
W(\omega_0)=\lim_{t \nearrow \infty}e^{ itH_{\omega_0}}e^{ -it( h_0
+\omega _0)\sigma_3}$ and set
$$R_+=W(\omega_0)^{-1}e^{i\theta(0)(P_+(\omega _0)- P_-(\omega
_0))}w_+  .$$ Notice that since $e^{ it\omega _0\sigma_3}$ is a
unitary matrix periodic in $t$ and $e^{i t\omega _0\sigma_3}R_+ $
describes circle in $L^2_x$, we have
$$\aligned & \lim_{t\to +\infty}
 \|e^{-itH _{\omega _0}} W(\omega_0)e^{i t\omega
_0\sigma_3}R_+ -e^{ -it( h_0 +\omega _0)\sigma_3}e^{i t\omega
_0\sigma_3}R_+\|_{H^1_x}=0.
\endaligned $$
Since $W(\omega _0)$ conjugates $H_{\omega _0}$ into $ \sigma
_3(h_0+\omega _0)$,  we get
$$e^{( it\omega
_0+i\theta (0))(P_+(\omega _0)- P_-(\omega _0))} e^{-itH _{\omega
_0}} w_+ =e^{-itH _{\omega _0}} W(\omega_0)e^{ it\omega
_0\sigma_3}R_+ .$$ Hence the last two limits  and the definition of
$R_+$ imply the limit
$$  \lim_{t\to +\infty} \left \|  e^{i\theta (t) \sigma _3}f (t) -
 e^{ -it h_0   \sigma_3}R_+  \right \|_{H^1_x}=0.$$
 Since $R(t)=\sum _j (z_j(t)\xi _j (\omega (t)) +\overline{z}_j(t)
 \sigma _1\xi _j (\omega (t)))+f(t)$ and since $\lim _{t\nearrow \infty}
 z_j(t)=0, $ it follows that
$$  \lim_{t\to +\infty} \left \|  e^{i\theta (t) \sigma _3}R (t) -
 e^{ -it h_0   \sigma_3}R_+  \right \|_{H^1_x}=0.\tag 1$$
Finally, by the same argument of \cite{BP2,BS,CM}   $\lim
_{t\nearrow \infty}\omega (t)=\omega _+$  for some $\omega _+$.
Hence expressing (1) in components we obtain Theorem 1.1.

\head \S 5 Review on Bloch functions for $h_0$ \endhead

By hypothesis for $\sigma (h_0)$   there exist $A_0^\pm
=0<A_1^{-}<A_1^{+ }<...<A_{ n_0 }^{-}<A_{ n_0 }^{+}<\infty$ such
that
$$\sigma (h_0)=\cup _{j=0}^{n_0-2}[A_{ j  }^{+},A_{ j+1
}^{-}]\cup [ A_{ n_0 }^{+}, \infty ).$$ The sets $ (A_j ^{-},A_j
^{+})$ for $j=1,...,n_0$ are nonempty spectral gaps. For $j\ge 0$ we set
$a_j^{\pm}= \sqrt{A_j^{\pm}}$. For  $j\le 0$,
$a_j^{\pm}=-a_{-j}^{\pm}$. For any $E  \in  {\Bbb C }\backslash
[0,\infty )$   there is a unique $ k\in {\Bbb C_+}$, called
quasimomentum such that $h_0u=Eu$ has a solution of the form $
\widetilde{\phi}_{\pm} (x,k)= e^{\pm ikx}\widetilde{m} _{\pm}(x,k)$
with $\widetilde{m}_{\pm}(x+1,k)\equiv \widetilde{m}_{\pm}(x,k)$ and
$\widetilde{m}_{\pm}(0,k)=1$. The correspondence $E\to k(E)$ is a
conformal map from $ {\Bbb C }\backslash [0,\infty )$ into a set
$$\Cal K=\Bbb C_+ \backslash ( \cup _{j=1}^{n_0} \left ([ n(j)\pi ,
n(j)\pi +ih _j ] \cup [ -n(j)\pi , -n(j)\pi +ih _j ] \right ) $$
with $n(j)$ strictly increasing for $j=1,.., n_0$ and $h_j>0$. We
set $\Cal T= \{ \pm \pi n(j): j=1,..., n_0\}$.
 The
inverse map $E=E(k)$ extends in  an even   map $\Bbb R\backslash \Cal T
\to \sigma (h_0)$. We have
$$E(  \pi n(j)  \pm 0 ):=\lim  _{\varepsilon \searrow 0}
 E(n(j)\pi \pm \varepsilon )=   \pi   A_j^{\pm}
.$$ $E(k)$ defined in $]\pi n(j), \pi n(j+1) [$ extends
  continuously  in $[ \pi n(j),
\pi n(j+1)]$. Same holds for $[\pi n(n_0), \infty )$. We consider
the
   solutions $\theta (x,k)$ and
$\varphi (x,k)$ of $h_0u=E(k)u$ which satisfy  the initial
conditions $ \varphi (0,k)=   \theta ^\prime (0,k) =0$  and $
\varphi ^\prime (0,k) =\theta (0,k)=1. $   Then $$\tilde \phi _\pm
(x,k)=\theta (x,k)+m^\pm (k) \varphi (x,k) \text{ with   } m^\pm
(k)=\frac{\varphi ^\prime (1,k) - \theta (1,k) }{2\varphi (1,k)} \pm
i \frac{\sin k}{\varphi (1,k)}.
$$
For real $k\in \Bbb R \backslash \Cal T$ we have $
\widetilde{\phi}_{\pm} (x,k)= \overline{\widetilde{\phi}_{\mp}
(x,k)} $. Set now
$$N^2(k)= \int_0^{1} \tilde \phi _+ (x,k)\tilde \phi _- (x,k) dx.$$
Then $N^2(k)>0$ for $k\in \Bbb R\backslash \Cal T$, with $N^2(k) $
bounded away from 0 near $\Cal T$, and $N^2(k)\neq 0$ for all $k\in
\Cal K$, see for instance \S 2 \cite{C2}. Then there is an $a>0$
such that  there is a well defined holomorphic square root $N(k)$
for $k\in \Cal K$ with $\Im k<a.$ Notice that the restriction $\Im
k<a $ is due to monodromy issues. We extend $N(k)\neq$ also for $\Im
k\ge a $, possibly as a discontinuous function. Then  we introduce
functions $m _\pm ^0(x,k) =\tilde m_{\pm }^{0} (x,k)/N(k) $  and
$\phi _{\pm }^{0} (x,k):=e^{\pm ikx}m_{\pm }^{0} (x,k)$. The
following  fact is well known.

\proclaim{Theorem 5.1} The   functions $\phi _{\pm }^{0}(x,k) $ are
analytic   in $k\in \Bbb R \backslash \Cal T$
 and for $ \Cal F _{h_0} {f} (k)=\int _\Bbb R    {\phi _+ ^{0}(x,k)}
f(x)dx$ we have:
$$\align & \int  _\Bbb R| f(y)|^2dy=  \int  _{\Bbb R}|
\Cal F _{h_0}f (k)|^2dk; \tag a
\\& f(x)=   \int _{\Bbb R}    \phi  _-^{0}
(x,k)  \Cal F _{h_0} f(k) dk; \tag b
\\& \Cal F _{h_0} \left ( h_0f  \right )(k)=
 E (k)  \Cal F _{h_0}f(k)    \text{ for $k\in \Bbb R \backslash \Cal T$}; \tag c \\&
e^{ith_0} (x,y)=\int   _{\Bbb R} e^{i(t E  (k) -(x-y)k) }   {m _-
 ^{0}(x,k)}{ m  _+  ^{0}(y,k)} dk .\tag d \endalign
$$

\endproclaim

In the sequel $\dot g (k)=\partial _k g(k)$. Korotyaev \cite{K}
proves:

\proclaim{Lemma  5.2} Consider    $E(k)$   for $k\in [ \pi n(j), \pi
n(j+1)]$. Then $ \dot E     (k)=0$ for $k=n(j)\pi , n(j+1)\pi  $ and
$\dot  E   (k)>0 $ in $]\pi n(j), \pi n(j+1) [$.   In $[ \pi n(j),
\pi n(j+1)] $ the equation $\ddot E     (k)=0$ admits exactly one
solution $k_j$. We have $k_j\in ] \pi n(j), \pi n(j+1) [$ and
$\dddot E
  (k_j)\neq 0$. We have $ \dot E   (k)>0 $ in
$]\pi n(n_0), \infty )$ and $\ddot E     (k)>c_0>0$  in $[\pi
n(n_0), \infty )$ for some $c_0>0$.
\endproclaim
Notice that $E(-k)=E(k)$.

\proclaim{Theorem  5.3}Set $e^{ith_0} (x,y)=\Cal U(t,x,y)+\Cal
V(t,x,y) $ with
$$\aligned & \Cal U(t,x,y)=\int   _{|k|\le \pi n(n_0)} e^{i(t E  (k) -(x-y)k) }   {m _-
 ^{0}(x,k)}{ m  _+  ^{0}(y,k)} dk\\&
\Cal V(t,x,y)=\int   _{|k|\ge \pi n(n_0)} e^{i(t E  (k) -(x-y)k) }
{m _-
 ^{0}(x,k)}{ m  _+  ^{0}(y,k)} dk.\endaligned
 $$
Then $|\Cal U(t,x,y) | \le C         \langle t\rangle    ^{-\frac
13}$ and $|\Cal V(t,x,y) | \le C         |t|  ^{-\frac 12} .
 $
\endproclaim
Assuming Lemma 5.4 below, the proof is  in  \cite{F1}. In \cite{C2}
  Theorem 5.3 is extended to the case of infinitely
many energy bands. In the case of finitely many energy bands we
have:

\proclaim{Lemma  5.4} For all $\alpha \in \Bbb N\cup \{ 0\}$ there
are constants $C_{\alpha}>0$ such that for $(x,y)\in \Bbb R^2$ and
for $k\neq \pm n(j)\pi$, $j=1,...,n_0$, we have
$$\aligned & \left | \partial _k^\alpha ( m ^0_+(x,k)  m ^0_-(y,k) -1 )\right | \le    C_\alpha \langle
k\rangle ^{-1 -\alpha }.
\endaligned$$
\endproclaim
{\it Proof (sketch).} Cases $\alpha = 0,1$ are stated and proved  in
\cite{C2} in the case of a generic periodic potential $P(x)$, i.e.
all gaps non empty. The case with finitely many bands is much
easier. In the finite bands case, there is an $M>0$ such that  the
$\alpha = 0 $ case is valid for all complex $k$ with $|k|\ge M$.
Here we are using the fact that $m ^0_+(x,k)  m ^0_-(y,k)$ is
holomorphic in $\Cal K \cup  \Cal K )^\ast \cup (\Bbb R \backslash
\Cal T),$ see \cite{C2}  and where $z\in \Cal K ^\ast$ if and only
if $\overline{z}\in \Cal K  $. The inequalities for $\alpha \ge 1$
hold by the Cauchy integral formula.

\bigskip In the
sequel we will set $w_0(k)=W[\phi _+^0 (x,k),\phi _-^0 (x,k) ].$

\proclaim{Lemma  5.5} For $k\neq \pm n(j)\pi$, $j=1,...,n_0$ we have
$\dot E(k)=-i w_0(k).$
\endproclaim
{\it Proof.} The formula $ \dot E(k)=\frac{2\sin (k)}{\varphi
(1,k)N^2(k)}$ is formula (3.1) \cite{F2}. By direct computation
$W[\widetilde{\phi }_+^0 (k),\widetilde{\phi} _-^0 (k)
]=i\frac{2\sin (k)}{\varphi (1,k) }$. So Lemma 5.5 follows by the
normalization in the definition of $\phi _\pm ^0 (x,k)$.

\head \S 6  Addition to $h_0$ of a small potential $q(x)\in C^\infty
_0(\Bbb R)$
\endhead
We consider the operator   $h=h_0+q(x)$ with $q(x)\in C^\infty
_0(\Bbb R)$     small. We consider solutions $f _{\pm } (x,k)=e^{\pm
ikx }m_{ \pm } (x,k)$ of $hu=E (k) u$ with

$$ \lim _{x\to +\infty } \frac {m_{ + } (x,k)}{ m _+^0(x,k) }=1 =
\lim _{x\to -\infty } \frac {m_{- } (x,k)}{ m _-^0(x,k) } . \tag
6.1$$

\proclaim{Lemma 6.1} Assume that $q(x)\in C^\infty _0(\Bbb R)$. Then
for any $k\in \Cal K$ the   equation

$$\aligned & m_{ \pm } (x,k)= m^0 _\pm  (x, k)-
\int _x^{\pm \infty }e^{\mp ik(x-t)}A(x,t,k) q(t) m_{ \pm } (t,k)
dt,\endaligned \tag 6.2
$$  $$\text{ with  }A(x,t,k):=\frac{ \phi _+  ^{0}(x ,k)\phi
_- ^{0}(t ,k)-\phi _-  ^{0}(x ,k)\phi _+  ^{0}(t ,k)   } {W[ \phi _+
^{0}(\cdot  ,k),\phi _- ^{0}(\cdot  ,k) ] },\tag 6.3 $$ has a unique
solution $m_{ \pm } (x,k)$ such that $f_{ \pm } (x,k)=e^{\pm ikx}m_{
\pm } (x,k)$ solves   $hu=E (k) u$ with the
 asymptotic property in (6.1). There is a constant $C_q$ such that
 for some $a>0$

$$  \aligned &|m_{ \pm } (x,k)-m^0 _\pm (x, k)|\le C_q
 \langle k \rangle ^{-1} {\max } ( 1,\mp x)  \\&
 \big |   \partial ^n  _km_{ \pm } (x,k)\big
|\le C_q  \langle k \rangle ^{-1} {\max } ^{n+1}( 1,\mp x)   \text{
for $n=1,2$  and $\Im k<a$}.
\endaligned \tag 6.4 $$ Furthermore we have continuous
 extensions for $k\in [ \pi n(j),
\pi n(j+1)]$, for  $k\in [\pi n(n_0),  \infty )
 $, for $k\in [ -\pi n(j+1),
-\pi n(j )]$ and for  $k\in (-\infty , -\pi n(n_0)]
 $ where the above estimates continue to be satisfied.
\endproclaim
Notice that estimates for $f_\pm (x,k)$ are   also in \cite{F2}. Set
$$D_{k}(x,t )=A(x,t,k) e^{-ik(x-t)}=\frac{ { m} _+ ^{0} (x ,k) {m} _- ^{0} (t ,k)-e^{-2ik(x-t)}
  {m} _- ^{0} (x ,k) {m} _+ ^{0} (t ,k) }
 {W[ {\phi} _+  ^{0}(x ,k),  {\phi} _-  ^{0}(x ,k) ]
} .$$ Then Lemma 6.1 follows from standard arguments, see Lemma 1
\cite{DT}, by the following lemma:

 \proclaim{Lemma 6.2} Let $t\ge x$. Then we have the
following for a fixed $C>0$:

{\item{(1)}}  $|D_{k}(x,t )|\le C \langle k \rangle ^{-1}$ for
$|k|\ge \pi n(n_0)+1$    ;
  {\item{(2)}}   $|D_{k}(x,t )|\le C \langle  x-t\rangle  $;
   {\item{(3)}}
$|\dot D_{k}(x,t )|\le C \langle x-t\rangle ^2 $  for  $\dot
D_{k}(x,t )=\partial _k D_{k}(x,t )$; {\item{(4)}} $|\ddot D_{k}(x,t
)|\le C \langle x-t\rangle ^3 $  for $\ddot D_{k}(x,t )= \partial
_k^2 D_{k}(x,t )$.

\endproclaim
(1) follows from Lemma 5.4 which yields a bound   $|m _\pm  ^{0} (x
,k)m _\mp  ^{0} (t ,k)|\le C$ for a fixed $C$ and $w_0(k)\approx
2ik$ for $|k|\nearrow \infty $.   $w_0(k)=0$ iff $k=k_0\pm 0$ for
$k_0\in \Cal T$.    By Lemmas 5.2 and 5.5 $  \dot w_0(k_0 \pm 0)
\neq 0  .
  $ So $|w_0(k)|\gtrsim \text{dist}(k,\Cal T)$. For fixed $t\ge x$
  consider   $x_1,t_1\in [0,2]$ with
$0\le t_1-x_1 \le t-x$ and  $t-x=t_1-x_1+  L$ with $L\in \Bbb N$.
Then we can set $   D_{k}(x,t )=D_{k}^{(1)}(x,t )+D_{k}^{(2)}(x,t )$
with

$$\aligned & D_{k}^{(1)}(x,t )=
\frac{  {m} _+ ^{0} (x _1,k) {m} _- ^{0} (t_1 ,k)-
 e^{ 2ik(t_1-x_1)} {m} _- ^{0} (x_1 ,k) {m} _+ ^{0} (t_1 ,k) } {w_0(k)
} \\&  D_{k}^{(2)}(x,t )= e^{ 2ik(t_1-x_1)}\frac{  1-  e^{ 2ikL}
 } {w_0(k)
}  {m} _- ^{0} (x_1 ,k) {m} _+ ^{0} (t_1 ,k).\endaligned
$$
$D_{k}^{(2)}(x,t )$ satisfies (1)--(4). Indeed Lemma 5.4 yields
bounds on ${m} _- ^{0} (x_1 ,k) {m} _+ ^{0} (t_1 ,k)$. Similarly,
one can bound $\frac{  1-  e^{ 2ikL}
 } {w_0(k)
} $ using $|w_0(k)|\gtrsim \text{dist}(k,\Cal T)$ and the fact that,
  for $k$ near $k_0\in  \Cal T$, then   $ w_0(k)$ for $k\ge k_0$
  resp.
  for $k\le k_0$ admits an analytic extension defined around $k_0$,
  see \cite{C2,F1,F2,K}. Finally we claim that $ |\partial _k^nD_{k}^{(1)}(x,t )|\le C$
  for $|n|\le 2$ and a fixed $C$. For $k$ large this is a
  consequence of $w_0(k)=i\dot E(k)\approx 2ik$ and the estimates on
  $E(k)$, see \cite{K}, and Lemma 5.4.   For $k_0\in \Cal T$ and $k_0^\pm = k_0\pm 0$

   $$\text{ $ e^{ik_0x }m _+ ^{0} (x
,k_0^+)\equiv e^{-ik_0x }m _- ^{0} (x   ,k_0^+)$ and $e^{ik_0x }m _+
^{0} (x  ,k_0^-)\equiv e^{-ik_0x }m _- ^{0} (x   ,k_0^-)$.}$$ Then
both numerator and denominator in the fraction defining
$D_{k}^{(1)}(x,t )$ are 0 for  $k=k_0\pm 0$ with $k_0\in \Cal T$.
Since, once again, these are functions admitting analytic extensions
on fixed neighborhoods  around $k_0$, one obtains the desired
estimates $ |\partial _k^nD_{k}^{(1)}(x,t )|\le C$.

\proclaim{Lemma 6.3} There exists  $q(x)\in C^\infty _0(\Bbb R)$ so
that for all $k\in \Cal T  $ with $k>0$  we have
$$ \align  & \int _{\Bbb R} q(x)| m
 _+^0(x,0  )|^2dx>0  \tag 1\\&
\mp \int _{\Bbb R} q(x)| m^0
 _+(x,k\pm 0  )|^2dx>0  \text{ for  all $k\in \Cal T$ with
 $k>0$}.\tag 2
  \endalign $$
 \endproclaim
{\it Proof.} $  m^0
 _+(x,k\pm 0  ) \not \equiv 0$ for all $k\in \Cal T$.
Since $\sigma (h_0)$ has finitely many bands, the potential $P(x)$
is analytic, and hence also the
 $m_\pm  ^0(x,k)$ are analytic in $x$.
(2) can fail if and only if   $| m^0
 _+(\cdot ,k_1+ 0  )|= | m^0
 _+(\cdot ,k_2- 0  )|$ for two   $k_1,k_2\ge 0 $  in $\Cal T$.
 The latter is the same of  $| \phi ^0
 _+(\cdot ,k_1+ 0  )|= | \phi ^0
 _+(\cdot ,k_2- 0  )|$.
Notice that $E(k_1+ 0 )\neq  E(k_2- 0 )$, so one of them is nonzero.
It is not restrictive to assume $E(k_2- 0 )>0$. Notice that $\phi ^0
 _+(x ,k_1+ 0  )$ and $\phi ^0
 _+(x ,k_2- 0  )$ are real valued for $x\in \Bbb R$. Indeed, for
 instance, by definition
$\overline{\widetilde{\phi} ^0
 _+(x ,k_1+ 0  )}= \widetilde{\phi} ^0
 _-(x ,k_1+ 0  )$, see (2.5) \cite{C2}, $w_0(k_1+ 0)=0$ implies  $\widetilde{\phi} ^0
 _+(x ,k_1+ 0  ) = C\widetilde{\phi} ^0
 _-(x ,k_1+ 0  )$ for a fixed $C$ and  by definition $\widetilde{\phi} ^0
 _+(0 ,k_1+ 0  ) =  \widetilde{\phi} ^0
 _-(0 ,k_1+ 0  )=1.$
 By $E(k_2- 0 )>0$ the function $\phi ^0
 _+(x ,k_2- 0  )$ has infinitely many zeros (if it has 1 it has infinitely
many since $\phi ^0
 _+(x ,k_2- 0  )=e^{ik_1x}m ^0
 _+(x ,k_2- 0  )$ with $m ^0
 _+(x ,k_2- 0  )$ periodic in $x$, if it had none $\phi ^0
 _+(x ,k_2- 0  )$ would correspond to $E(k_2)=0$).
 We choose two points $a<b$ with  $\phi ^0
 _+(a ,k_2- 0  )=\phi ^0
 _+(b ,k_2- 0  )=0$  and $\phi ^0
 _+(x ,k_2- 0  )\neq 0$ for $a<x<b$. $| \phi ^0
 _+(\cdot ,k_1+ 0  )|= | \phi ^0
 _+(\cdot ,k_2- 0  )|$ implies the same statement   for
$\phi ^0
 _+(x ,k_1+ 0  )$.  But then they are both ground states for
 the Dirichlet  problem in $[a,b]$ for the operator $h$. Then
$\phi ^0
 _+(x ,k_1+ 0  )=C \phi ^0
 _+(x ,k_2- 0  )$ for fixed $C$, but this is impossible by
$E(k_1+ 0 )\neq  E(k_2- 0 )$.

\bigskip

 We pick  $q(x)\in C^\infty _0(\Bbb R)$ small as in Lemma 6.3.
By \cite{FK} we have that $\sigma (h)=\sigma (h_0)$, and in
particular $h$ has no eigenvalues. Furthermore there are no
resonances, that is there are not nonzero solutions of any of the
equations  $hu=A_j^\pm u$  for $j>0$ with $u\in L^\infty$. Combining
the asymptotics in Lemma 6.1 and the analysis of the thresholds in
\cite{FK} we conclude:

 \proclaim{Lemma 6.4} For the Wronskian  we have $W[f_+(k), f_-(k)]\neq 0$
  for all $k$ and $W[f_+(k),
 f_-(k)]\approx 2ik$ for $|k|\to \infty .$
\endproclaim
\bigskip
  There are equalities
$$ \aligned &f_{ \mp } (x,k)= \frac {R_{\pm } (k) }{T (k)} f_{ \pm }
(x,k) +\frac {1 }{T (k)} f_{ \pm } (x,-k)\\& \overline{f_{ \pm }
(x,k)}=f_{ \pm } (x,-k)\, , \, \overline{T(k)}=T(-k) \, , \,
\overline{R_{\pm } (k)}=R_{\pm } (-k)
\\&  |T(k)|^2+ |R_{\pm } (k) |^2=1\, , \quad T(k)\overline{ R_{\pm }}
(k)+R_{\mp} (k) \overline{T(k)}=0 \endaligned \tag 6.6
$$
with   $T(k)$ and $R_{\pm }(k)$ smooth functions for $k\in \Bbb R
\backslash \Cal T
 $  with smooth extensions in $  [ \pi n(j),
\pi n(j+1)]$,   $  [\pi n(n_0),  \infty )
 $,   $  [ -\pi n(j+1),
-\pi n(j )]$ and   $   (-\infty , -\pi n(n_0)]
 $ and with $T (k)=0$  for all $k\in \Cal T$   and $T(0)=0$.   For $q(x)$ small
 $$\big |  [T(k)-1] ^{(m)}\big | + \big |
  R_{\pm }^{(m)}(k)\big |\le  {\delta }/{\langle k\rangle
 } \text{ for $|m|\le 2$}.\tag 6.7$$

 Next define

$$ \psi  (x ,k)= \left \{ \matrix   \frac{1}{\sqrt{2\pi  }}T(k) f_{ + } (x,k) \quad \text{for}
\quad k \ge 0 \\  \frac{1}{ \sqrt{2\pi }}T(-k) f_{ - } (x,-k) \quad
\text{for} \quad k < 0 \endmatrix \right   .\tag 6.8
$$
Then by standard arguments (partially repeated in the proof of Lemma
8.12) we have the following version of Theorem 5.1 for $h$.
\proclaim{Lemma 6.5} For $ \Cal F _h {f} (k):=\int _\Bbb R \psi (x
,k) f(x)dx$  we have

$$\align & \int  _\Bbb R| f(y)|^2dy=  \int  _{\Bbb R}|\Cal F _h f (k)|^2dk \tag a
\\& f(x)=   \int _{\Bbb R} \overline{\psi  (x ,k)} \Cal F _h {f} (k) dk \tag b
\\& \Cal F _h\left( {h f  }\right ) (k)= E (k)  \Cal F _h {f}(k)     \tag c \\&
e^{ith } (x,y)=\int   _{\Bbb R} e^{i t E  (k)    }   \overline{\psi
(x ,k)}  {\psi  (y ,k)}  dk .\tag d \endalign
$$
\endproclaim
We have:
 \proclaim{Lemma 6.6 (dispersive estimates)  }
Set $e^{ith } (x,y)=\Cal U(t,x,y)+\Cal V(t,x,y) $ with
$$\aligned & \Cal U(t,x,y)=\int   _{|k|\le \pi n(n_0)} e^{i t E  (k)
   }   \overline{\psi
(x ,k)}  {\psi  (y ,k)}  dk\\& \Cal V(t,x,y)=\int   _{|k|\ge \pi
n(n_0)} e^{i t E  (k)    }   \overline{\psi (x ,k)}  {\psi  (y ,k)}
dk.\endaligned
 $$
Then $|\Cal U(t,x,y) | \le C         \langle t\rangle    ^{-\frac
13}$ and $|\Cal V(t,x,y) | \le C         |t|  ^{-\frac 12}
 $ for a fixed $C>0$.
\endproclaim
For $\Cal U(t,x,y)$ the argument is the same of Theorem 5.3. Turning
to $\Cal V(t,x,y)$, the proof follows along the lines of Lemmas 2.4
and 2.5 in \cite{Wd2} using \cite{C2}, so we give a short sketch
only. Let $ \varsigma (k)$ the characteristic function of the set
$|k|\ge \pi n(n_0).$ Set
 $\phi _\pm (t,x,y)=\int   _{\Bbb R ^{\pm}} e^{i t E  (k) }
\overline{\psi (x ,k)} {\psi  (y ,k)} \varsigma (k) dk$. Split $\phi
_\pm =\phi _\pm ^{low}+\phi _\pm ^{high}$ distinguishing large and
small $k$ in the integral. To fix ideas let us consider only
 $\phi  (t,x,y)=2\pi \phi _+ ^{high}(t,x,y)=2\pi \int   _{\Bbb R^{+}}
  e^{i t E  (k)    }   \overline{\psi
(x ,k)} {\psi  (y ,k)} \varphi  (k) dk$ with $\varphi   $ a  smooth
increasing function equal to 0 near $\pi n(n_0)$ and equal to 1 near
$+\infty$.  It is not restrictive to assume $x<0$ and $y>0$ since
with different signs we can use (6.6) as in \cite{Wd2}. Then $\phi =
\phi _1+\phi _2$ with the $\phi _j$ defined below.   We have $\phi
_1=\phi _1^{(0)}+\phi _1^{(1)}$ with
$$\aligned & \phi _1 ^{(0)}(t,x,y)=\int   _{\Bbb R^{+}} e^{i t E  (k) +ik(x-y)   }
{m _-^{0} (x ,k)} {m  _+^0(y ,k)} \varphi  (k)  dk\\& \phi _1
^{(1)}(t,x,y)=\int _{\Bbb R^{+}} e^{i t E  (k) +ik(x-y)   } {m
_-^{0} (x ,k)} {m  _+^0(y ,k)}  \varphi  (k)
(|T(k)|^2-1)dk.\endaligned
$$
The proof of $| \phi _1 ^{(0)}(t,x,y)|\le Ct^{-\frac 12}$ follows
from the arguments in the proof of Lemma 4.3 in \cite{C2} which in
part we repeat now to prove $| \phi _1 ^{(1)}(t,x,y)|\le Ct^{-\frac
12}$.  We extend $E  (k)$ outside the support of $\varphi$ so that
is a strictly convex smooth even function on $\Bbb R$. Let $p(k)=E
(k) +ikt^{-1}(x-y)$. Let $k_0$ be the unique solution of $\dot p (k)
=0$. For $\chi (t)$ a cutoff supported near $t=0$, we
 insert the partition of unity  $\chi (k-k_0)+(1-\chi (k-k_0))$.
 Correspondingly set $\phi _1 ^{(1)}=\phi _1 ^{(1,1)}+\phi _1
 ^{(1,2)}$
with
$$\Phi _1 ^{(1,1)}=\int _{\Bbb R^{+}} e^{i t E  (k) +ik(x-y)   } {m
_-^{0} (x ,k)} {m  _+^0(y ,k)}  \varphi  (k) (|T(k)|^2-1)\chi
(k-k_0)dk$$ and $\phi _1 ^{(1,2)}$ defined similarly but with
$(1-\chi (k-k_0)).$ Then by stationary phase $ \left | \phi _1
^{(1,1)}(t,x,y)\right | \le Ct^{-\frac 12}.$  Set $q^2/2=p (k) -p
(k_0)$. Then
$$\aligned &\Phi _1
^{(1,2)}(t,x,y)=e^{i  tp (k_0) }
  \int   _{\Bbb R}
e^{itq^2 }  \rho (q) dq \text{  where}\\& \rho (q)={m _-^{0} (x
,k(q))} {m _+^0(y ,k(q))} \varphi (k(q))(1-\chi
(k(q)-k_0))(|T(k(q))|^2-1)  {dk}/{dq}.
\endaligned$$ Then
$   | \phi _1 ^{(1,2)}(t,x,y)  | \le C  t ^{-\frac{1}{2}} \| \rho
(q)\| _{H^{ 1}}$.   Inequality $\|  \rho (q)\| _{H^{ 1}}<\infty$
follows by Lemma 5.4, by (6.7) and by $\| (1-\chi (k(q)-k_0))
dk/dq\| _{H^{ 1}}<\infty $, Lemma 4.9 \cite{C2}. We have
$$\aligned & \phi _2  (t,x,y)=\int   _{\Bbb R^{+}}
 e^{i t E  (k) +ik(x-y)   }
|T(k)|^2g_{x,y}(k) \varphi  (k)  dk \text{ with}\\& g_{x,y}(k) ={m
_-  (x ,k)} {m _+ (y ,k)}- {m _-^{0} (x ,k)} {m _+^0(y ,k)}
.\endaligned
$$
Then also $   | \phi _2(t,x,y)  | \le C  t ^{-\frac{1}{2}}$ by (6.4)
and proceeding as above.

\head \S 7 Wave operators and partial   diagonalization for $  H
_\omega $
\endhead

\noindent We   write $  H _\omega =  \sigma_3(h+\omega )+  B^\ast
(\omega ) A $ with $A(x  )=\langle x \rangle ^{-\tau}$ with $\tau
>3/2$ and
$B^\ast (x,\omega )$ a    $C^2$ function   in $(x,\omega )$ with
values in the space of $2\times 2$ real valued matrices. For any
$\omega  $ in some compact set $K$  there is a  constant  $c (K)>0$
and $\alpha
>0$ such that
$   \big | e^{\alpha |x|}  B^\ast (x,\omega ) \big | \le c_m(K)
 $ $  \forall   x\in \Bbb R. $  We have:

 \proclaim{Proposition 7.1} Assume that $H_\omega$ does not have
 resonances at $\pm \omega$  and $\sigma _e(H_\omega )$ does not
 contain eigenvalues.
Then there are isomorphisms inverses of each other
 $W(\omega )\colon  L^2_x \to L^2_c(  H_{\omega }) $ and
$Z(\omega )\colon L^2_c(  H_{\omega })  \to  L^2_x
  $,
defined as follows:
  \noindent for $u\in L^2_x$, and $v$ such that $
  \sigma _3v\in L^2_c(  H_{\omega }
  )  $,
$$\aligned & \langle W u,v\rangle =
\langle  u,v\rangle  +\lim _{\varepsilon \to 0^+ } \frac 1 {2\pi i }
\int _{-\infty}^{+\infty} \langle A R_{\sigma_3(h+\omega )}(\lambda
+i\varepsilon ) u,B R_{H_\omega ^\ast }(\lambda +i\varepsilon
)v\rangle d\lambda ;
\endaligned$$
for $u\in L^2_c(  H_{\omega }) $,  $v\in  L^2_x $,
$$\aligned & \langle Z u,v\rangle =
\langle  u,v\rangle  +\lim _{\varepsilon \to 0^+ } \frac 1{2\pi i }
\int _{-\infty}^{+\infty} \langle A R_{H_\omega }(\lambda
+i\varepsilon ) u,B R_{\sigma_3(h+\omega )}(\lambda +i\varepsilon
)v\rangle d\lambda .
\endaligned$$ Then $ P_c(  H _\omega
)  H _\omega =W \sigma_3(h+\omega )Z.$
  $\| W(\omega )\| _{B(L^2_x,L^2_c(  H_{\omega }   )
)}$ and   $\| Z(\omega )\| _{B(L^2_c(  H_{\omega }   )  , L^2_x) }$
are uniformly locally bounded in $\omega $.
\endproclaim
We need to show that  there is a fixed  $c>0$ such that $\forall \,
\epsilon \neq 0$
$$\align  & \int \|  \langle x \rangle ^{-\tau} R_{\sigma_3(h+\omega )}(\lambda +i\varepsilon) u\|^2_{L^2_x} d\lambda \le
c \| u\| ^2_{L^2_x}  \text{ for all $u\in L^2_x$}\tag 1\\& \int \|
B R_{\sigma_3(h+\omega )}(\lambda +i\varepsilon) u\|^2_{L^2_x}
d\lambda \le c \| u\| ^2_2\text{ for all $u\in L^2_x$}\tag 2\\&
\int \|  B R_{H_\omega ^\ast }(\lambda +i\varepsilon ) u\|
^2_{L^2_x} d\lambda \le c \| u\|^2_{L^2_x} \text{ for all $u\in
L^2(H _\omega ^\ast ):=\sigma _3L^2(H _\omega   )$} \tag 3
\\& \int \|  \langle x \rangle ^{-\tau} R_{H_\omega }(\lambda
+i\varepsilon)  u\|^2_{L^2_x} d\lambda \le c \| u\|^2_{L^2_x}
\text{ for all $ u \in L^2_c(H _\omega ) $} . \tag 4
\endalign
$$
Let us first prove (1), (2). They are consequences of (5) for $\tau
>3/2$:
$$\int \|  \langle x \rangle ^{-\tau} R_{h}(\lambda +i\varepsilon) u\|^2_{L^2_x} d\lambda \le c \|
u\|^2_{L^2_x}  \text{ for all $u\in L^2_x$}.\tag 5$$
 By (5.3) in
Theorem 5.1 \cite{K},  (5) will follow from   $$ \| \langle x
\rangle ^{-\tau}R_{h } (z ) \langle x \rangle ^{-\tau} \|
_{B(L^{2}_x, L^2_x)}<C \text{  for all $z$ with $0<|\Im z |  $}.\tag
6$$ Observe that for $  \Im k\ge 0$, $\pm (x-y)\ge 0$ and setting
$w(k)=W[f_{+}( k), f_{-}(  k)]$,

$$   \langle x \rangle ^{-\tau} R_{h}( E(k),x,y)
 \langle y \rangle ^{-\tau}=       {e^{ik|x-y|}
\langle x \rangle ^{-\tau} m_{\pm}(x,k)        {m_{\mp}(y,k)}\langle
y \rangle ^{-\tau}} w^{-1}(k).$$ By Lemma 6.1, $|m_{\pm}(x,k)
m_{\mp}(y,k)|\lesssim \langle x \rangle \langle y \rangle $. By
Lemma 6.4, $ |w(k)|\gtrsim \langle k \rangle .$    Hence
 $$\langle x \rangle ^{-\tau} |R_{h}( z ,x,y)|
 \langle y \rangle ^{-\tau}\le C{\langle z \rangle}^{-\frac{1}{2}}
\langle x \rangle ^{1-\tau  }
 \langle y \rangle ^{1-\tau} \tag 7$$ for all $z$ with $\Im z\ge 0$
  resp.
 $\Im z\le 0$ (for $\Im z=0$ there are
 two different continuations). Then for all $z$ with $\Im z\ge 0$
 resp.
 $\Im z\le 0$

$$ \| \langle x \rangle ^{-\tau} R_{h}( z ,x,y)
 \langle y \rangle ^{-\tau} \| _{L^2 _{x,y}}
  \le  {C}{\langle z \rangle}^{-\frac{1}{2}} . \tag 8$$
 (8) implies (6). We consider now (3) and (4),
 and for definiteness we focus on (4).  We have for $A=\langle x \rangle
 ^{-\tau}$
$$AR_{H_\omega}(z)u=
 ( 1+ AR_{\sigma_3(h+\omega )}(z)B^\ast   )^{-1}
  AR_{\sigma_3(h+\omega )}(z)u.\tag 9
$$
(8) implies $\| AR_{\sigma_3(h+\omega )}(z)B^\ast \|
_{B(L^2,L^2)}\lesssim \langle z \rangle ^{-1/2}.$ By Fredholm theory
there is a bounded 0 measure set $X\subset \Bbb R$ such that in
$\Bbb R \backslash X$   the following limits  exist
$$ \lim _{\epsilon \searrow 0}R_{H_\omega }(\lambda +i\varepsilon )
 = R_{H_\omega }^\pm  (\lambda )  \text{ in $C^0_{loc}(\Bbb
R \backslash X, B(L^{2, \tau }_x ,L^{2, -\tau }_x )) $}.$$ A point
$\lambda $ belongs to $X$ exactly if
  $\ker ( 1+ AR_{\sigma_3(h+\omega )}^{\alpha}(\lambda )B^\ast   ) \neq
  0$ in $L ^{2 }$  for $\alpha =+$ or $-$. The points in
  $X\backslash \sigma _e(H_\omega )$ are  eigenvalues of $H_\omega$.
  This follows from (9). Furthermore, the hypothesis that
$\sigma _e(H_\omega )$ contains no eigenvalues of $H_\omega$ implies
that $X\backslash \sigma _e(H_\omega )$ is exactly the set of
eigenvalues of $H_\omega$. The exponential decay of $H_\omega -
\sigma _3(h_0+\omega )$ implies, by standard arguments,  that
$H_\omega $ has finitely many eigenvalues. Therefore for $z$ close
to
  an eigenvalue of  $H_\omega $ and for $u\in L^2_c(H_\omega)$   we have that
  $\| AR_{H_\omega}(z)u\| _{L^2_x}\le c
  \|  R_{H_\omega}(z)u\| _{L^2_x}\le c' \| u\| _{L^2_x}$
  for fixed constants. We will show below that $X\cap \sigma _e(H_\omega )$
is empty. Notice that this yields (4).   To see this let $N \subset
\Bbb R$ be the neighborhood of $\sigma _e(\Cal H_\omega ) $ formed
by the points with distance $< \delta $ from $\sigma _e(\Cal
H_\omega ) $, with $\delta >0$ a small number. We  split the
integral in (4) in two integrals with domain $\Bbb R\backslash N$
and $N$. The integral with domain $\Bbb R\backslash N$ is bounded,
since for $\lambda \in \Bbb R\backslash N$ we have $\| \langle x
\rangle ^{-\tau} R_{  H_\omega }(\lambda +i\varepsilon) u\| _{
2}\le \| R_{  H_\omega }(\lambda +i\varepsilon ) u\| _{ 2}\le C  \|
u\| _{ 2}$ for a fixed $C$  and $\Bbb R\backslash N$ is a bounded
set. To bound the integral with domain $N$ we use formula (9)  above
where $H_0=\sigma_3(h+\omega )$
$$\langle x \rangle
 ^{-\tau}R_{ H_\omega}(\lambda +i\varepsilon )u =
 ( 1+ \langle x \rangle
 ^{-\tau}R_{  H_0}(\lambda +i\varepsilon )B^\ast   )^{-1}
  \langle x \rangle
 ^{-\tau}R_{  H_0}(\lambda +i\varepsilon ) u.
$$
Then for any $u\in L ^2$
$$\align  &  \| \langle x \rangle
 ^{-\tau}R_{  H_\omega}(\lambda +i\varepsilon )u \| _{L^2_
\lambda (N, L ^2_x)} \le \\& \| ( 1+ \langle x \rangle
 ^{-\tau}R_{  H_0}(\lambda
+i\varepsilon )B^\ast   )^{-1} \| _{L^\infty_ \lambda (N,B(L ^2_x, L
^2_x))} \| \langle x \rangle
 ^{-\tau}R_{
H_0}(\lambda +i\varepsilon ) u \| _{L^2_ \lambda (N, L ^2_x)}.\tag
10
\endalign  $$
The first factor in rhs(10)  is uniformly bounded by the hypothesis
that $
  X\cap \sigma _e(H_\omega )$
is empty, while the
 second factor in
rhs(10) is bounded for any $u\in \ell ^2$ by (1).  This yields (4),
with the proof of (3)   similar.

Now we return to the set $X$. Our aim is tho show that $X$ is formed
exactly by the eigenvalues of $H_\omega $. It is enough to show that
$X\cap \sigma _e(\Cal H _\omega )$ is empty. We will proceed in two
steps. We will first show that by    (H7)    $X$ does not contain
any of these extremes thresholds of  $\sigma_3(h+\omega )$,
  Lemma 7.2. We
  will then show that the points of $X$ in the interior of the
  spectral bands are necessarily eigenvalues of $H_\omega$. Since
  such "embedded" eigenvalues do not exist by (H7), then we
  can conclude that
$X\cap \sigma _e(\Cal H _\omega )$ is empty.

  \proclaim{Lemma 7.2}
$X$ does not contain   extremes   of the spectral bands of
  $\sigma_3(h+\omega )$.
\endproclaim
  {\it Proof.} Suppose
the claim   is wrong, and that   $\lambda \in X$ is an extremum of
one of the spectral bands of
  $\sigma_3(h+\omega )$, and pick $w\neq 0$  in $L ^{2  }$
  such that
  $$( 1+ AR_{\sigma_3(h+\omega )} (\lambda )B^\ast   )w=0.$$
Then $( 1+ B^\ast AR_{\sigma_3(h+\omega )} (\lambda )  )B^\ast w=0
$. Then $\psi = R_{\sigma_3(h+\omega )} (\lambda ) B^\ast w$ is
$\psi \neq 0$ and by standard arguments is
 a
nontrivial distributional solution of $(H_\omega -\lambda )u=0$. We
claim that $$\psi \in L^\infty_x .\tag 10$$ If (10) is true, we get
contradiction with hypothesis (H7). The proof of (10) reduces at
showing that for a rapidly decreasing bounded function $g(y)$ we
have
$$\| f\| _{L^\infty_x}<\infty \text{ for }f(x)=
\int _{\Bbb R}R_{ h } (x,y,\lambda )g(y)dy .$$ This is a consequence
of Lemmas 6.1 and 6.4 and, for $w(k)=W[f_+(k),f_-(k)]$, of
$$f(x)=f_{+}(x,k)\int _{-\infty}^{x}f_{-}(y,k)g(y)dy/w(k)+
f_{-}(x,k)\int ^{ \infty}_{x}f_{-}(y,k)g(y)dy/w(k).$$ This proves
  Lemma 7.2.

   \proclaim {Lemma 7.3} For any $\lambda _0$ in the interior of
   $\sigma _e(H_\omega )$ there are an open  interval
   $I\subset \sigma _e(H_\omega )$ with $\lambda _0\in I$
   and a constant $C$ such that
   $\| AR(\lambda \pm i\varepsilon )u\| _{L^2_\lambda (I)L^2_x}
   \le C \|  u\| _{ L^2_x}$ for any $u$.
   \endproclaim
{\it Proof.} Let   $AB^\ast  _1 =AB^\ast -\sigma _3 q$. Then
$H_\omega = \sigma _3(h_0+\omega )+AB^\ast  _1 $. The points
$\lambda _0  \in \sigma (\sigma _3(h_0+\omega )) $
  where $\| AR(\lambda  \pm i \varepsilon  )u\| _{L^2_\lambda (I)L^2_x}$
   is unbounded  for $I$ a small neighborhood of $\lambda _0$
   are such that if we set
$$AR_{H_\omega}^\pm (\lambda )u=
 ( 1+ AR_{\sigma_3(h_0+\omega )}^\pm (\lambda )B^\ast  _1 )^{-1}
  AR_{\sigma_3(h_0+\omega )}(\lambda )u, \tag 11
$$
  we have for one of the two signs $$\ker ( 1+
AR_{\sigma_3(h_0+\omega )}^\pm (\lambda _0)B^\ast _1 )\neq 0.$$ By a
standard argument, see \S 2 \cite{CPV}, if we assume   Lemma 7.4
below we conclude that such point $\lambda _0$   is an eigenvalue of
$ H_\omega$. Since we are excluding eigenvalues inside the spectral
bands of $\sigma _e( H_\omega )$ we are done with the proof Lemma
7.3.

\bigskip Lemmas 7.2 and 7.3 imply that  $X$ coincides with
$\sigma _d(H_\omega )\cap \Bbb R$. The above arguments moreover
imply that the following limits exist
$$ \lim _{\epsilon \searrow 0}R_{H_\omega }(\lambda \pm i\varepsilon )
 = R_{H_\omega }^\pm  (\lambda )
 \text{ in $C^0 ( \sigma _e(H_\omega ),
  B(L^{2, \tau }_x ,L^{2, -\tau }_x )) $}. \tag 12$$
\bigskip
\proclaim{Lemma 7.4} Fix $E_0$ in the interior of a spectral band of
$h_0$. Suppose that there is a function $\psi (x)$ such that  $\psi
\in L^{2,s}_x$ for all $s>0$ and $\Cal F _{h_0}(\psi )(\pm k_0)=0$,
where $E_0=E(k_0)$. Then $R^{\pm}_{h_0}(E_0)\psi \in L^{2 }_x$.
\endproclaim
{\it Proof.} Consider a cutoff $\chi (E)\in C^\infty _0(\Bbb R)$
with $\chi (E)=1$ for $E$ close to $E_0$ and $\chi (E)=0$ if
$|E-E_0|\ge \varepsilon _0>0$ for $\varepsilon _0$ a small fixed
number. Then
$$ R^{\pm}_{h_0}(E_0)\psi = R^{\pm}_{h_0}(E_0)\chi (h_0) \psi +
R^{\pm}_{h_0}(E_0)(1-\chi (h_0)) \psi .$$ The second term in the
right is in $L^2_x.$ With the notation in \S 5 we have
$$\Cal F _{h_0}(\chi (h_0) \psi )(  k ) =
\int _\Bbb R e^{ikx}m _+ ^0(x,k)\chi (E(k))  \psi (x) dx.$$ Notice
that  $m _+ ^0(x,k)\chi (E(k))$ is the symbol of a smoothing
pseudodifferential operator. This implies that  $\partial ^n _k\Cal
F _{h_0}(\chi (h_0)\psi ) \in L^2_k $ for all $n$ by Lemma 5.4.
Since $k_0\not \in \Cal T$, $\dot E (k_0)\neq 0$ by Lemma 5.2. Then
also
$$\partial ^n _k\left [ \frac{\Cal F _{h_0}(  \psi )(  k )
\chi (E(k)) }{E(k)-E_0} \right ] \in C^0 \cap L^1_k \text{ for all
$n$}.
$$ We have
$$R^{\pm}_{h_0}(E_0)\chi (h_0) \psi (x)=
\int _\Bbb R e^{-ikx}m _- ^0(x,k) \frac{\Cal F _{h_0}(\psi )(  k
)}{E(k)-E_0} \chi (E(k))  dk.$$ Integrating by parts we conclude $ x
^N R^{\pm}_{h_0}(E_0)\chi (h_0) \psi (x) \in L^\infty_x $ for all
$N$.

\bigskip
By Proposition 7.1 and by the spectral Theorem we conclude:

\proclaim{Lemma 7.5} For any  $u\in L^2_x$ we have
$$\aligned & P_c( H _{\omega} )u=\lim _{M \to \infty }\lim _{\epsilon \to 0^+}
 \frac 1{2\pi i}  \int _{\sigma _e(H_\omega )}
 \chi _{[-M,M]}(\lambda )\left [ R_{H_\omega }(\lambda +i\epsilon
)- R_{H_\omega }(\lambda -i\epsilon )
 \right ]  ud\lambda   .
\endaligned
$$
\endproclaim
Finally, we obtain the limiting absorption principle:

\proclaim{Lemma 7.6} For any  $u\in L^{2,\tau }_x$ with $\tau
>3/2$
we have
$$\aligned & P_c( H _{\omega} )u=\lim _{M \to \infty }
 \frac 1{2\pi i}  \int _{\sigma _e(H_\omega )}
 \chi _{[-M,M]}(\lambda )\left [ R_{H_\omega }^+(\lambda
)- R_{H_\omega }^-(\lambda   )
 \right ]  ud\lambda   .
\endaligned
$$
\endproclaim
{\it Proof.} For $u\in L^{2,\tau }_x$ and for fixed $M$ the
$\epsilon \to 0^+$ limit in Lemma 7.5 converges  in $L^{2,-\tau }_x$
to the integral in Lemma 7.6.

  \head \S 8 Plane waves  for $H_\omega $ \endhead

We recall $H_\omega =\sigma _3(h_0+\omega ) +V(x)$ with:

{\item{$\bullet $} }  $V(x)$ a real entries square 2 matrix s.t.
$\sigma _1V(x)=-V(x)\sigma _1$ and $\sigma _3V(x)= {^tV}(x)\sigma
_3$; {\item{$\bullet $} }  $|V(x)|\le C e^{-\sqrt{2\omega }|x|}$
 ; {\item{$\bullet $} } $H_\omega
$ has no eigenvalues in $\sigma _e(H_\omega )$ and the points and
the thresholds of the spectral bands of $\sigma _e(H_\omega )$ are
not resonances. This last statement means that for any $k\in \Cal
T$, and for $\lambda =E(k\pm 0)+\omega )$, as well as for $\lambda
=-E(k\pm 0)-\omega  $,   if $g\in L^\infty_x$ satisfies
$H_\omega g=\lambda  g$, then $g= 0$.

We set
$$\eta (x)=\int _x^\infty |V(t)|dt, \quad \gamma (x)=\int _x^\infty
\langle t \rangle |V(t)|dt.$$ We have $\sigma _e(H_\omega )= (\omega
+\sigma (h_0 ))\cup (-\omega -\sigma (h_0 ))$. Because of the
symmetries of $H_\omega$ we look at Jost functions $ {F_\pm (x,k)}$
close to $\omega +\sigma (h_0 )$.

 For   $\lambda \not \in \sigma _e(H_\omega ) $ let
  $k\in \Bbb C$  with $\Im k \ge 0$   such that   $\lambda =E(k)+\omega $
  (in the sequel, for $\lambda$ and $k$ in the same sentence, we will have
  always the relation $\lambda =E(k)+\omega $).
  We consider now   solutions  of $  H
_\omega u=\lambda u$ of the following
 form, where $ A(x,t,k) $ is given by (6.4):
$$\aligned   {F}_{  \pm } (x,k)=&  \phi _\pm  ^{0}(x,k)\overrightarrow{e}_1 - \int _x^{\pm \infty
} A(x,t,k)
\text{diag }(1,0) V(t)  {F}_{ \pm } (t,k) dt-\\
-&  \int _\Bbb RR_{h_0 }(-2\omega -E (k), x,t)
 \text{diag
}(0,1) V(t)  {F}_{  \pm } (t,k) dt. \endaligned \tag 8.1
$$

 \proclaim{Lemma  8.1} Assume (H7) and (H9).
 Then there is a small $\delta
>0$ and a finite set $S\subset \overline{\Bbb C} _+$
such that for $k\not \in S$ with $0\le \Im k\le \delta $,  (8.1) has
for any choice of sign a unique solution satisfying the estimates
listed below  for $M_{  \pm } (x,k):=e^{\mp ikx}{F}_{ \pm } (x,k)$.
These solutions solve $ H_\omega u=\lambda u$ with $\lambda =\omega
+E(k)$ and with the asymptotic property $ {F}_{ \pm } (x,k)= {\phi
}_{  \pm }^{0} (x,k)\overrightarrow{e}_1+o(e^{\pm ikx})$ for $x\to
\pm \infty .$
    For any fixed $0<a<\sqrt{2\omega }$
    and any $ \varepsilon >0$ there
 is
  $C_\varepsilon$    such that
$\forall \, x\in \Bbb R$ and
  $\forall \, k \in \Bbb C$ with $0\le \Im k\le \delta $ and    $\text{dist}(k,S)>\varepsilon >0$
  we have:

$$  \align &  | M_{  \pm }
(x,k)- m_{  \pm }^0 (x,k)\overrightarrow{e} _1 |\le C_\varepsilon
e^{-a \max (1,\mp x )} \langle k \rangle ^{-1} (1+\max (1,\mp x ) )
; \tag 1
\\&   \big |  {\partial
_k } [ M_{\pm  } (x,k)-m_{  \pm }^0 (x,k)\overrightarrow{e} _1 ]\big
|\le C_\varepsilon   {\max }^2(1,\mp x ) ;\tag 2 \\ & \big |
{\partial ^2_k } [ M_{\pm } (x,k)-m_{  \pm }^0
(x,k)\overrightarrow{e} _1 ]\big |\le C_\varepsilon   {\max
}^3(1,\mp x ) .\tag 3 \endalign
$$
\endproclaim
  {\it Proof.}  We sketch the proof for the $+$ case
   and drop the index. Using the
notation in \S 6,
$$\aligned  {M} (x,k)=&   {m} ^{0}_+(x,k)\overrightarrow{e}_1 - \int _x^{  \infty
}  D_k (  x,t )
\text{diag }(1,0) V(t)  {M}  (t,k) dt-\\
& - \int _\Bbb RR_{h_0 }(-2\omega -E (k), x,t) e^{
-ik(x-t)}\text{diag }(0,1) V(t)
 {M}  (t,k) dt. \endaligned \tag 8.2
$$
We have
 $ \big | R_{h_0 }(-2\omega -E (k), x,t) e^{
-ik(x-t)} \big |  \le C    e^{-  \sqrt{ 2\omega }  |x-t| }\langle
k\rangle  ^{-1} .$ Let  for $V=\{  V_{\ell, j}\}$ $$\aligned &
^tM=(M_1,M_2)\, , \quad M_1 ^{(1)}(x,k):= - \int _x^{ \infty } D_k (
x,t )V_{12}(t) M_2 (t,k) dt  \, ,  \\& M_1 ^{(0 )}(x,k):= m
_{+}^{0}(x,k)  - \int _x^{ \infty } D_k ( x,t )V_{11}(t) M_1 (t,k)
dt.\endaligned $$ Then $M_1 (x,k)= M_1 ^{(0 )}(x,k) + M_1
^{(1)}(x,k)$. By standard arguments:

\proclaim{Lemma 8.2} For given $ M_2\in L^\infty _x$   for some
fixed $C=C(V )$ we have:

$$ \align  & |M_1 ^{(0 )}(x,k)-m_+ ^{0}(x,k)|\le C
  {\max } (1,- x )  \gamma (x) \langle k\rangle  ^{-1}  ;\tag 1
\\&  |M_1 ^{(1)}(x,k) |\le
C  {\max } (1,- x )  \gamma (x) \langle k\rangle  ^{-1} \|  M_2
(\cdot ,k)\| _{L^\infty_x }    . \tag 2
\endalign
$$
\endproclaim

\bigskip
{\it Continuation of proof of Lemma 8.1.} We have $M_1 ^{(1)}(x,k)=L
(k)M_2 (x,k)$  with $L(k)$ a linear operator such that for a fixed
$C$
$$\aligned & |L(k)M_2( x)     |\le    C     \langle k   \rangle ^{-1} {\max } (1,- x )
    \| M_2
   (\cdot ,k)\| _{L^\infty_x }  .\endaligned  \tag 3$$
   Eliminate $M_1$ from     (8.2) to get the following system
for $K(k)$ defined below

$$\aligned  &(1+K(k)) M_2 (x,k)=- \int _\Bbb RR_{h_0 }(-2\omega -E (k), x,t) e^{
-ik(x-t)}  V_{21}(t)  M_1 ^{(0)}(t ,k)
  dt \endaligned\tag 8.3$$ $$ K(k)g(x) =  \int _\Bbb RR_{h_0 }(-2\omega -E (k), x,t) e^{
-ik(x-t)} \left ( V_{21}(t)L (k)[g] (t)+ V_{22}(t) g(t) \right )
  dt.    $$

\proclaim{Lemma 8.3} Let $0\le  \Im k  < \delta < \sqrt{2\omega} -a
$ for some $a>0$.  Then $K(k)$ maps $L^\infty_x$ in
$e^{-a|x|}W^{1,\infty }_x$. There is a finite set $S\subset
\overline{\Bbb C} _+$ such that for $k$ as above with $ k\not \in S$
equation (8.3) admits a unique solution $M_2(x,k)\in L^\infty_x $.
Furthermore, for any $\varepsilon >0$ there is a $C_\varepsilon$
such that if $\text{dist}(k,S)>\varepsilon  $ then

$$ \left  | e^{a|x|} M_2 (x,k)\right  |  \le  {C _\varepsilon}\la k
\ra ^{-1} .\tag 1$$
\endproclaim
{\it Proof. } The fact that $K(k):L^\infty_x \to e^{-a|x|}W^{1,\infty
}_x
 $  if $0\le  \Im k  < \delta < \sqrt{2\omega} -a
$ is an elementary computation. For $|k|\to \infty $
 we have $ \| K(k)\| _{B(L^\infty_x , L^\infty_x )}   \to 0$
 so we can   solve (8.3) obtaining (1).
  $ K(k):L^\infty_x \to L^\infty_x$ is a compact
operator. For $|k|\lesssim 1$ by  the  Fredholm alternative $S$ will
be formed by the $k$ with   $ \ker (1+K(k))\neq 0$. The set $S$ is
necessarily discrete and contained in $|k|\lesssim 1$. We start by
considering an  $k\in S\cap \Bbb R.$ If we have such a $k$ then, we
have a nonzero solution of the equation obtained from (8.3)
replacing $M_1 ^{(0)}(t ,k)$ with 0. Going backwards in the above
argument, we obtain the existence of a nontrivial $ g \in L^\infty_x$ satisfying
$$\aligned {g}  (x )=&    - \int _x^{  \infty
} A(x,t,k)
\text{diag }(1,0) V(t)  {g}  (t ) dt-\\
-&  \int _\Bbb RR_{h_0 }(-2\omega -E (k), x,t) \text{diag }(0,1)
V(t) {g}  (t ) dt.\endaligned \tag 8.4
$$
Such $g\in L^\infty_x$ is a distributional solution of
$H_\omega g =\lambda g$ with $\lambda= E(k) +\omega  $. By
hypothesis we are assuming such solutions do not exist for
 $\lambda =E(k\pm 0 )+\omega $ with $k\in \Cal T$,
that is for   $\lambda$ a  threshold in $\sigma _e(H_\omega )$. By
continuity, $ \ker (1+K(k))= 0$ near such thresholds. Since the
points in $\Cal T$ are the only possible accumulation points of $S$,
we conclude that $S$ is finite.

\bigskip
\bigskip
{\it Continuation of proof of Lemma 8.1.} Lemmas 8.2-3 imply (1)  in
Lemma 8.1.  To prove (2) we differentiate (8.2) in $k$ obtaining for
the dot representing $\partial _k$, and for $E_k(x,t)=R_{h_0
}(-2\omega -E (k), x,t)$,

$$\aligned  & \dot {M} (x,k)= \overrightarrow{e}_1 \dot m_+^0(x,k)-
\\& - \left ( \int _x^{  \infty
}  D_k (x,t) \text{diag }(1,0) +\int _\Bbb R E_k(x,t)\text{diag
}(0,1) \right )   V(t)
 \dot {M}  (t,k) dt \\& - \left ( \int _x^{  \infty
}  \dot D_k (x,t) \text{diag }(1,0) +\int _\Bbb R \dot
E_k(x,t)\text{diag }(0,1) \right )   V(t)
  {M}  (t,k) dt . \endaligned \tag 4
$$
We have $  |\partial _k^j E_k(x,t) |\le C e^{ -\alpha \langle k
\rangle |x-t|   }  $ for $j=0,1$ for fixed $C>0$ and $\alpha >0$. We
write   $\dot {M} _1(x,k) = \dot {M} _1^{(0)}( x,k)+\dot {M}
_1^{(1)}( x,k)$, with $\dot {M} _1^{(1)}( x,k)= L(k)\dot {M} _2(
x,k) $,  with

$$\aligned & \dot {M} _1^{(0)}( x,k)\overrightarrow{e}_1=
\overrightarrow{e}_1 \dot m_+^0(x,k)+\\&
  +\left (\text{ third line of (4)}\right )-\overrightarrow{e}_1
  \int _x^{ \infty } D_k (  x,t )V_{12}(t) \dot M_2 (t,k)
dt  .\endaligned$$ Then  $$\aligned & |\dot {M} _1^{(0)}( x,k)|\le C
{\max }^2 (1,-x) \text{ and }\\&     | L(k)\dot {M} _2( x,k)|\le
C\langle k \rangle ^{-1} {\max }  (1,-x) \| \dot {M} _2( \cdot ,k)\|
_{L^\infty_x} .\endaligned $$ We obtain an analogue of system (8.3), with
$ {M} _1^{(0)}( x,k)$ replaced by $\dot {M} _1^{(0)}( x,k)$. Then we
conclude $   | e^{a|x|} \dot M_2 (x,k) | \le \frac{C }{1+|k|} . $
Repeating the argument in a similar way, we get claim (3) in Lemma
8.1.

\bigskip

If $S\cap \Bbb R=\emptyset$, by taking $\delta$ small enough we can
neglect the set $S$. However in Lemma 5.4 \cite{C1} the proof that
$S\cap \Bbb R=\emptyset$ is wrong.   In fact the discussion in
\cite{KS} allows for the existence of $S$. So let us assume now
$S\cap \Bbb R\neq \emptyset$. For $k$ near $S\cap \Bbb R$  we
consider the system

$$ (1+R _{\sigma _3(
h _0+\omega )}^{+ }(\lambda )V) {\Phi}_{ \pm } (\cdot ,k)=\phi ^{0}
_\pm
 (x,k)
 \overrightarrow{e}_1  . \tag 8.5$$
 \proclaim{Lemma 8.4 } Let $\lambda =E(k)+\omega $ with  $ E(k)\in \sigma (h_0) $ with
$k\neq \pm \pi n(j)$ for all $j$.

{\item{(1)}}For   any choice of signs, system (8.5) admits exactly
one solution in $  L^\infty_x$. {\item{(2)}} For $k\not \in
S$ we have $ \Phi _{ \pm } (x,k)= c_{\pm} (k) F _{ \pm } (x,k)$ for
some constants $c_{\pm} (k)$. {\item{(3)}}Fix a decomposition
$V=B^\ast A$ with $A$ and $B^\ast$ $\in C^2$ and exponentially
decreasing. Then $k\to A\Phi _{ \pm } (x,k) \in L^2_x $ is a real
analytic map for $k$ near $S$. {\item{(4)}} $c_\pm (k)$ are real
analytic functions in $k$ for $k$ near $S$.

\endproclaim
{\it Proof.} For definiteness pick $+$. For $\lambda =E+\omega$ we
claim that $R _{\sigma _3( h _0+\omega )}^{+ }(\lambda )V $
$$  \text{ is a   $L^\infty_x\to
  L^\infty_x$ compact operator.}\tag 5$$
 We have
$$R _{\sigma _3(
h +\omega )}^{+ }(\lambda )V=\sigma _3\left [  \matrix R _{h _0}^{+
}(E )V_{11} & R _{h _0}^{+ }(E )V_{12}\\
R _{h _0}^{+ }(-E-2\omega )V_{21} & R _{h _0}^{+ }(-E-2\omega
)V_{22}
\endmatrix\right ]  .$$
So the proof of (5) reduces at proving that  the operators $R _{h
_0}^{+ }(-E-2\omega )\nu(x)$ and $R _{h _0}^{+ }(E ) \nu (x)$ are
compact from $L^\infty_x$ into itself for $\nu (x)$ an
exponentially decreasing $C^2$ scalar function. These operators are
defined in $L^\infty_x$ with values in $W^{1,\infty }_x$.
Since $R _{h _0}^{+ }(-E-2\omega )\nu (x)$ maps $L^\infty_x$
in $e^{-a|x|}L^\infty_x$ for a fixed $a=a(\omega _0)>0$,
it is easy to conclude that $R _{h _0}^{+ }(-E-2\omega )\nu (x)$
satisfies (5). We have for $w_0(k)=[\phi _+^0(k),\phi _-^0(k)]$, $R
_{h _0}^{+ }(E )[\nu g] (x)=$
$$\aligned & =\phi _+^{0}(x,k) C_- (g)-  \int _{x}^{ \infty} \left
[ \phi _+^{0}(x,k)  \phi _-^{0}(y,k) - \phi _+^{0}(y,k)  \phi
_-^{0}(x,k)\right ] \nu(y) g(y) dy/w_0(k) \\& =\phi _-^{0}(x,k) C_+
(g)+ \int ^{x}_{-\infty} \left [ \phi _+^{0}(x,k)  \phi _-^{0}(y,k)
- \phi _+^{0}(y,k)  \phi _-^{0}(x,k)\right ] \nu(y) g(y) dy/w_0(k)
\endaligned$$
with $  C_\pm   (g)=\int _\Bbb R \phi _\pm ^{0}(y,k)\nu(y)g(y)
dy/w_0(k) .$   Given a sequence $\{ g_n \} $ with $\| g_n \|
_{\infty}\le 1$ we want to show that a subsequence of $\{ R _{h
_0}^{+ }(E )[\nu g_n] \}$ converges in $L^\infty_x$. It is
not restrictive to assume that $\{ R _{h _0}^{+ }(E )[\nu g_n] \}$
converges in $L^\infty_{loc} (\Bbb R)$, and  $\{ C_+ (g_n) \} $ and
$\{ C_- (g_n) \} $ converge in $\Bbb C$. It is elementary conclude
by the above formulas on $R _{h _0}^{+ }(E )[\nu g] (x) $, that $\{
R _{h _0}^{+ }(E )[\nu g_n] \}$ converges in $L^\infty_x$.

Having established that (5) is compact, by Fredholm theory we know
that (8.5) has a unique solution unless there is a nonzero  solution
$g\in   L^\infty_x$ of
$$ (1+R _{\sigma _3(
h _0 +\omega )}^{+ }(\lambda ) V )g=0  . \tag 6$$ Notice that $^tg
=(g_1,g_2)$

$$H_\omega g=\lambda g \text{ in the sense of  distributions}  \tag 7$$
and that $|g_2(x)|\le C e^{- \sqrt{2\omega }|x|}$ for some $C>0$. By
(6) we obtain
$$\langle (1+R _{\sigma _3(
h _0+\omega )}^{+ }(\lambda )V ) g, \sigma _3Vg\rangle =0.$$ Then
$$ \langle \delta (h _{0} -E(k)) (V_{11}g_1+V_{12}g_2 ),
V_{11}g_1+V_{12}g_2 \rangle  =0.\tag 8$$ Notice that
$V_{11}g_1+V_{12}g_2 \in L^{2,s}_x$ for any $s\in \Bbb R$. For $\psi
\in L^{2,s}_x$ with $s>1/2$ we have
$$\langle \delta (h _{0} -E(k)) \psi ,
\psi \rangle  =|\dot E (k))| ^{-1}\left ( | \Cal F _{ h_0}\psi
(k)|^2+| \Cal F _{ h_0}\psi (-k)|^2\right ).$$ So if we set   $\psi
=V_{11}g_1+V_{12}g_2 $, then $\Cal F _{ h_0}\psi (\pm k)=0$. We can
apply Lemma 7.4 and conclude that
$$g_1=R _{
h  _{0}}^{+ }(E(k))\psi\in L^{2 } (\Bbb R) . $$ So $g\in L^{2 }_x$
and (7) imply that $g$ is an eigenfunction and that $\lambda$ is an
eigenvalue. Since we are excluding this possibility, we conclude
that there are non nonzero solutions of (6) and that Claim (1) is
correct.

  From (8.5) we get for $w_0(k)=W[\phi ^{0} _+(k),\phi ^{0}
_-(k)]$
$$\aligned & \Phi _{ \pm } (x,k)=\phi ^{0}
_\pm
 (x,k)
 \overrightarrow{e}_1 \\& +\phi ^{0}
_+
 (x,k)
  \int _{-\infty}^{x}\phi ^{0}
_-
 (y,k) \text {diag} (1,0)  V (y)\Phi _{ \pm } (y,k) dy/w_0(k)\\& -\phi ^{0}
_-
 (x,k)
  \int^{ \infty}_{x}\phi ^{0}
_+
 (y,k) \text {diag} (1,0)  V (y)\Phi _{ \pm } (y,k) dy/w_0(k) \\
 & - \int _\Bbb RR_{h_0 }(\omega -E (k), x,y) \text{diag }(0,1) V(y)
\Phi _{ \pm } (y,k) dy.
\endaligned \tag 8.6$$
Then, for $$\aligned & c_+ (k)=1+\int _{\Bbb R} \phi ^{0} _-
 (y,k)  (V _{11}(y)\Phi _{1 +} (y,k)
 +V _{12}(y)\Phi _{2 + } (y,k)) dy/w_0(k)\\& c_- (k)=1-
 \int _{\Bbb R} \phi ^{0} _+
 (y,k)  (V _{11}(y)\Phi _{1 -} (y,k)
 +V _{12}(y)\Phi _{2 - } (y,k)) dy/w_0(k)  \endaligned \tag 9$$
 and for $k\not \in S$ we have for $g_\pm (x):=
 \Phi _{ \pm } (x,k) - c_\pm (k)F _{ \pm } (x,k)$,
$$\text{$g
\in L^\infty_x$ with }\lim _{x\to \pm \infty }g_\pm  (x )=0
\text{ and } H_\omega g_\pm  =\lambda g_\pm  .$$ Since $k\not \in S$
this yields $g_\pm \equiv 0$.

To prove claims (3) and (4) we write (8.5) as
$$ (1+AR _{\sigma _3(
h _0+\omega )}^{+ }(\lambda )B^\ast) A{\Phi}_{ \pm } (\cdot
,k)=A\phi ^{0} _\pm
 (x,k)
 \overrightarrow{e}_1  .  $$
Since $AR _{\sigma _3( h _0+\omega )}^{+ }(\lambda )B^\ast \in B
(L^2 _x,L^2 _x)  $ is real analytic in $k$, and similarly $k\to
A\phi ^{0} _\pm
 (x,k)\in L^2 _x$ is real
analytic in $k$, we obtain claims (3) and (4)
 by the chain rule.
\bigskip

\bigskip

Lemma 8.1 yields Jost functions $F(x,k)$ for $H_\omega $ for energy
$\lambda$ close to $ \omega +\sigma (h)$ and for $k\not \in S$. By
$\sigma _1 H _\omega =-   H _\omega \sigma _1 $ we conclude that
$\sigma _1F(x,k)$ are Jost functions for $H_\omega $ for energy
$\lambda$ close to $- \omega -\sigma (h)$. From $\sigma _3  H
_\omega  =  H _\omega ^\ast \sigma _3   $ we get that $\sigma
_3F(x,k)$ are Jost functions for $H_\omega ^\ast $ for energy
$\lambda$ close to $ \omega +\sigma (h)$ and $ \sigma _3\sigma
_1F(x,k)$ are Jost functions for $H_\omega ^\ast $ for energy
$\lambda$ close to $- \omega -\sigma (h)$.

\proclaim{Lemma 8.5} For   $k\in \Bbb R\backslash S$ we have
$\overline{F_{\pm}(x,k)} =F_{\pm}(x,-k)$ and for   $k\neq 0$ we have
$$ {F}_{  \mp } (x,k)=
\frac{1}{ T  (k  )}\overline{ {F}_{  \pm } (x,k)} +\frac{R_{ \pm}
(k)}{T  (k  )}  {F}_{  \pm } (x,k) \tag 1 $$ where $T  (k  )$ and
$R_\pm (k  )$ are defined by the above formula.
\endproclaim
{\it Proof.} $\overline{F_{\pm}(x,k)} =F_{\pm}(x,-k)$ follows  by
the fact that $V(x)$ has real entries and by uniqueness in Lemma
8.1. We claim that the triples of functions in formula (1) in the
statement are linearly dependent. Notice that this implies
immediately the statement, since one can see near $+\infty$ by Lemma
8.1 that the functions on the right hand side of (1) are linearly
independent. For definiteness we will prove linear dependence of $
{F}_{ +} (x,k) $,  $ {F}_{  - } (x,k)$ and $ \overline{{F}_{ - }
(x,k)}$. If  we assume they are linearly independent, we can find a
nonzero linear combination $g=\alpha F_++\beta F_-+\gamma
\overline{F_-}$ with $g(x)\in L^\infty_x$ and  $g(x)\to 0$ for $x\to
+\infty$. Then $g(x)$ satisfies (8.4) and by the fact that $k\not
\in S$ we conclude $g(x)=0$. So our assumption is absurd and the
three functions are linearly dependent.

 \proclaim{Lemma 8.6} All the
Wronskians below are constant. We have for any $k\in \Bbb
R\backslash S $:
$$  \align   & W[\overline{ {F}_{ \pm } (k)},    {F}_{  \pm } (k)]
=W[ \phi  _{\mp}^{0}(\cdot , k  ), \phi _{\pm}^{0} (\cdot ,k )] ,
\tag 1\\&
 T(k  )= \frac {W[ \phi  _{\mp}^{0}(\cdot , k  ), \phi
_{\pm}^{0} (\cdot ,k )]} {W[ F_{  \mp } (k),  F_{ \pm } (k)] } \, ,
\quad R_{  \pm } (k ) =  \frac {W[  {F_{ \mp } ( k )},
 \overline{F_{ \pm }} (k)] } {  W[ F_{ \pm
} (k), F_{  \mp } ( k ) ] }  \tag 2
\\&
  \overline{T  (k  )}=T  (-k  ) \, , \,
\overline{R_{  \pm } (k  )}=R_{  \pm } (-k  ) , \tag 3
\\&|T  (k  )|^2+ |R_{ \pm } (k ) |^2=1\, , \quad T (k)  \overline{R_{ \pm }
(k  )}+R_{  \mp} (k  ) \overline{ T  (k  )}=0.\tag 4\endalign
$$
\endproclaim
{\it Proof.} For the fact that the Wronskians are constant see Lemma
5.8 \cite{KS}. The rest of the proof is the same of Lemma 5.6
\cite{C1}.

\proclaim{Lemma 8.7} For $k\in \Bbb R$ we have
$$W[ \Phi _{  + } (k), \Phi _{ - } (k)] =c_+(k)
W[ \phi _{+}^{0}(\cdot , k  ), \phi _{-}^{0} (\cdot ,k )] = c_-(k)
W[ \phi _{+}^{0}(\cdot , k  ), \phi _{-}^{0} (\cdot ,k )].\tag 1$$
In particular $c_+(k)=c_-(k)= T(k).$
\endproclaim
{\it Proof.} We have the asymptotic behaviors, with $\partial
_xo(1)=o(1)$,
$$\aligned & \Phi _{  + } (x,k)=  c_+(k)\phi _{+}^{0} (x ,k )
\overrightarrow{e}_1+o(1) \text{ for } x\to +\infty\\& \Phi _{  - }
(x,k)=  \phi _{-}^{0} (x ,k )\overrightarrow{e}_1 + C_1(k)\phi
_{+}^{0} (x ,k ) \overrightarrow{e}_1+o(1)\text{ for } x\to +\infty
\\& \Phi _{  - } (x,k)= c_-(k)\phi _{-}^{0} (x ,k )
\overrightarrow{e}_1 +o(1)\text{ for } x\to -\infty\\& \Phi _{  + }
(x,k)=  \phi _{+}^{0} (x ,k )\overrightarrow{e}_1 + C_2(k)\phi
_{-}^{0} (x ,k ) \overrightarrow{e}_1+o(1)\text{ for } x\to -\infty
\endaligned $$ with $c_\pm (k)$ defined in Lemma 8.4 and for some
constants $C_j(k)$ with $j=1,2.$ In particular, for later use,
$C_1(k)$  is smooth for $k$ near $S$ and given by
$$\aligned & C_1(k)=
  \int _{\Bbb R}\phi ^{0}
_-
 (y,k)   ( V_{11} (y)\Phi _{ 1-} (y,k) +  V_{12} (y)\Phi _{ 2-} (y,k)
 )dy/w_0(k).
\endaligned \tag 2$$
(1) follows by using these asymptotic expansions and the fact that
the Wronskian $W[ \Phi _{ + } (k), \Phi _{ - } (k)]$ is a constant.
We have
$$\aligned & T(k  )= \frac {W[ \phi  _{-}^{0} , \phi
_{+}^{0}  ]} {W[ F_{ - }  ,  F_{ + }  ] } =c_+c_-\frac {W[ \phi
_{-}^{0} , \phi _{+}^{0}  ]} {W[ \Phi  _{ - }  ,  \Phi _{ + }  ] }
=\frac{c_+c_-}{c_\pm} =c_\pm .
\endaligned $$

 \proclaim{Lemma 8.8} $T (k)$ and $R_{ \pm }
(k)$   are in $C^1(\Bbb R)$ and there is $C
>0$ such that for $n\le 1$
 $$\big |  {d^n}\left [ T  (k)-1
 \right ] /{dk^n}\big | +  \big | {d^n R_{  \pm } (k )}/{dk^n}
  \big |   \le  {C  }/{\langle k\rangle }.
$$
\endproclaim
{\it Proof.} Notice that for $k\not \in T $    we have $W[ \phi
_{-}^{0}(k) , \phi _{+}^{0} (k) ] \neq 0$. So  (4) Lemma 8.6  and
the non resonance hypothesis  at the thresholds imply $W[ F_{ -
} , F_{ + } ] \neq 0$ in $\Bbb R \backslash S$. There $T  (k)$
 and $R_{  \pm } (k )$ are $C^1$ with the desired asymptotic
 estimates by Lemma 8.1. Near $S$ we use $T(k)=c_\pm (k)$ and (3)
 Lemma 8.4 to conclude that $T\in  C^1(\Bbb R)$. Near $S$
 $$\aligned & R_+(k)= \frac{W[F_-,\overline{F}_+]}{W[F_+, {F}_-]}=
\frac{|c_+|^2}{c_+ ^2} \frac{W[\Phi _-,\overline{\Phi }_+]}{W[\Phi
_+, {\Phi }_-]} =C_1(k) \frac{ \overline{c}_+ ^2}{c_+ ^2}.
\endaligned $$
Here $C_1(k)$ is in (2) Lemma 8.7 and is smooth. If $k_0\in S$ then
there is some $m \ge $ such that $c_+^{(m)}(k_0)\neq 0$, then also
$\overline{c}_+/c_+$ is smooth near $k_0$. So $R_+\in  C^1(\Bbb R)$.
The argument for $R_-$ is similar.

\bigskip
We consider the following system :
$$\aligned  {G} _{\pm} (x,k )=&    \phi _\pm ^{0}(x,k_1) \overrightarrow{e}_2  - \int _x^{ \pm  \infty
} A(x,t,k)
\text{diag }(1,0) V(t)  {G} _{\pm} (t,k )dt-\\
-&  \int _x^{\pm \infty} A(x,t,k_1)\text{diag }(0,1) V(t)  {g}
_{\pm} (t,k ) dt \endaligned \tag 8.7
$$
where $E(k_1)=-E(k)-2\omega $.

 \proclaim{Lemma 8.9} For $\lambda =\omega +E(k)$ with $0\le \Im k < \delta
 $
 system (8.7) admits exactly one solution $G_{\pm }(x,k)$ which satisfies $H_\omega G_{\pm }=\lambda G_{\pm
 }$, $G_{\pm }(x,k)$ is real for $k\in \Bbb R  $. There is
 $C $ such that
$$  \align &  | {G}_{  \pm }
(x,k)-  \phi _\pm ^{0}(x,k_1)  \overrightarrow{e}_2 |\le C
   \langle k \rangle ^{-1}   \left | \int _x ^{\pm \infty}\langle t \rangle |V(t)| dt \right |
   e^{\mp x\Re \sqrt{k^2+2\omega } } ;
\\ & \big |
\partial _x {G}_{  \pm } - \partial _x\phi _\pm ^{0}(x,k_1)
\overrightarrow{e}_2 |\le C    \left | \int _x ^{\pm \infty}\langle
t \rangle |V(t)| dt\right |  e^{\mp x\Re \sqrt{k^2+2\omega } }
  .   \endalign
$$
We have $w(k)=[ {G_{- }( k)},  G_{+ }( k)] $ is a continuous
function with $w(k)= 0$ if and only if $k\in S$. We have $w(k)
=2\sqrt{k^2+2\omega } (1+o(1))$ for $|k|\to \infty .$

\noindent We have $[ {F_{\pm }(x,k)},  G_{\pm }(x,k)]=0$ for any
$k\not \in S$. For $\Im k>0$ and  $g\in L^\infty (b,+\infty )$, for
some $b\in \Bbb R$, with $H_\omega g=\lambda g$, then $g (x)=\mu
F_{+ }(x,k)+ \nu G_{+ }(x,k)$ for constants $\mu$ and $\nu$. If
$g\in L^\infty (-\infty ,b)$ solves $H_\omega g=\lambda g$, then
$g(x)=\mu F_{- }(x,k)+ \nu G_{- }(x,k)$
\endproclaim
{\it Proof.} The fact that there is a unique $G_{\pm }(x,k)$
satisfying (8.6) and the estimates  follow  from the argument in
Lemma 1 \cite{DT}. The asymptotic expansion of $w(k)=[ {G_{- }( k)},
G_{+ }( k)] $ for $|k|\to \infty  $ follows from the inequalities.
$k\in S$ exactly when $G_\pm (x,k)\in L^\infty _x$. This is
the case exactly for $ [ {G_{- }( k)},  G_{+ }( k)] =0 $. Identity
  $[ {F_{\pm } }, G_{\pm } ]=0$ follows by using the fact
that the Wronskian is constant and by the asymptotic properties of
the functions. Consider $h\in L^\infty (b,+\infty )$ with $H_\omega
h=E h$ with $\Im k>0$ and $k\not \in S$. By standard arguments there
is a $b\le b_1$ such that in $[b_1,+\infty )$ there are two other
solutions $a(x,k)$ and $b(x,k)$ of $H_\omega u=E u$, both unbounded
and such that $a(x,k)$, $b(x,k)$, $F_+(x,k)$ and $G_+(x,k)$ form a
fundamental set of solutions. Then in $[b_1,+\infty )$ we have
$h(x)=\mu F_{+ }(x,k)+ \nu G_{+ }(x,k)$ and by unique continuation
this holds  in $(b ,+\infty )$.

\proclaim{Lemma 8.10} Let $\lambda =\omega +E(k)$ with $0<\Im
k<\delta  $  and $k\not \in S$.  The  resolvent $R_{H_\omega }
(\lambda )$ has integral kernel given by $R_{H_\omega } (x,y,\lambda
)=R_1(x,y )+R_2(x,y )$ with for $x<y$
$$\aligned &  R_{1 } (x,y ) =-\frac{F_{ - } (x,k)  { ^t(\sigma _3 F_{ + } (y,k)} )}
{W[ F_{ + } (\cdot ,k),    F_{ - } (\cdot ,k) ] } \text{ and }R_{2 }
(x,y ) =-\frac{G_{ - } (x,k)  { ^t(\sigma _3 G_{ + } (y,k)} )} {W[
G_{ + } (\cdot ,k),    G_{ - } (\cdot ,k) ] }
\endaligned $$
and for $x>y$ $$\aligned &  R_{1 } (x,y ) =-\frac{F_{ + } (x,k)  {
^t(\sigma _3 F_{ - } (y,k)} )} {W[ F_{ + } (\cdot ,k),    F_{ - }
(\cdot ,k) ] } \text{ and }R_{2 } (x,y ) =-\frac{G_{ + } (x,k) {
^t(\sigma _3 G_{ - } (y,k)} )} {W[ G_{ + } (\cdot ,k),   G_{ - }
(\cdot ,k) ] }.
\endaligned $$
We have $ R_{H_\omega }
(x,y,\overline{\lambda})=\overline{R_{H_\omega } ( x,y ,\lambda )} $
and $ R_{H_\omega } (x,y,-{\lambda })=-\sigma _1R_{H_\omega } (x,y,
{\lambda })\sigma _1.$

\endproclaim
{\it Proof.} The last two equalities follow from the fact that
$H_\omega $ is real and $\sigma _1 H_\omega =-H_\omega \sigma _1 .$
Let   $R(x,y)=R_1(x,y)+R_2(x,y)$. To
 show that $R (x,y)$ is the kernel of the resolvent it is enough to
 show that for any fixed $x$ and for any fixed $\lambda$ the following
 equalities hold:
 $$ \aligned & A(x,\lambda):=
 \partial _x R (x,x^-)-\partial _x R (x,x^+)=-\sigma
 _3  \\&
  B(x,\lambda):=R (x,x^-)-  R (x,x^+)=0.\endaligned \tag 1$$
Now, for $ h(x)= F_+,F_-,  G_+,G_- $, using the information on the
Wronskians, we get
$$ ^th (x) A(x,\lambda )-{^th '(x)}B(x,\lambda )=- {^th} (x)
\sigma _3.\tag 2$$ If now, for $(x,k)$ fixed,  $(h(x),h'(x))$
 span $\Bbb C^2\times \Bbb C^2$ for $ h(x)= F_+,F_-,  G_+,G_- $,
 then  (2) yields (1).
 For $|k|\gg 1$ this is the case. Since $A(x,\lambda )$ and
 $B(x,\lambda )$ depend analytically on $k$, (1) holds for all $k
 \not \in S$.

We next define:
 \proclaim{Definition 8.11(Plane Waves)} For $k\not \in S$ set $
\Psi (x ,k):=\frac{1}{\sqrt{2\pi} } T  (k  )  {F}_{ + } (x,k )$ for
$k\ge 0$,   $ \Psi (x ,k):=\frac{1}{\sqrt{2\pi} } T  (-k )
 {F}_{  - } (x,-k)$
for $k< 0$. Near $S$ set  $ \Psi (x ,k)=\frac{1}{\sqrt{2\pi}
}{\Phi}_{ + } (x,k )$  for $k>0$ and $ \Psi (x
,k)=\frac{1}{\sqrt{2\pi} }{\Phi}_{ - } (x,-k )$  for $k<0$.
\endproclaim
The interpretation of the $\Psi (x ,k)$ is justified by the
following lemma:

\proclaim{Lemma 8.12} The operator $P_c(H_{
 \omega}  )$ has, for $^tA$ the transpose of $A$, kernel
$$ P_c( H_{
 \omega} )(x,y)= \int _\Bbb R \left [ { \overline{\Psi  (x ,k)}} \, \, \,  {^t\left
( \sigma _3  {\Psi    (y ,k)}\right ) }
  +   \sigma _1 { \overline{\Psi  (x ,k)}}\, \,  \,  {^t\left (  \sigma _3\sigma
_1  {\Psi    (y ,k)}\right ) } \right ] dk ,$$   In particular we
have $ (e^{itH_\omega}P_c(H_{
 \omega} ))(x,y)=$
$$   \int _\Bbb R \left [ e^{it(E(k)+\omega )}{ \overline{\Psi  (x ,k)}} \, \, \,  {^t\left
( \sigma _3  {\Psi    (y ,k)}\right ) }
  +  e^{-it(E(k)+\omega )} \sigma _1 { \overline{\Psi  (x ,k)}}\, \,  \,  {^t\left (  \sigma _3\sigma
_1  {\Psi    (y ,k)}\right ) } \right ] dk .$$
\endproclaim

{\it Proof.} Recall $ \langle f ,
  g\rangle =\int _{\Bbb R}{^tf(x)} {g(x)} dx.$ For $f\in L^2_c  (H  _{\omega }) \cap
  \Cal S(\Bbb R)$
    we have
$$\aligned & \langle  f ,g\rangle =
\lim _{\epsilon \to 0^+} \frac{1}{2\pi i}\int  _{\sigma _e(H_\omega
)} \langle \left [ R _{H  _{\omega }}(\lambda +i\epsilon )- R_{H
_{\omega }}(\lambda -i\epsilon ) \right ] f,g \rangle d\lambda  \\&
=\frac{1}{2\pi i}\int _{\sigma _e(H_\omega )} \langle H(\lambda )f,g
\rangle d\lambda \text{ with $H(\lambda )f (x)=\int _{\Bbb R }
H(\lambda ,x,y)f (y)dy,$}
\endaligned $$
$$ H(\lambda ,x,y)= R_{H  _{\omega }}(\lambda +i0,x,y )-
   R_{H  _{\omega }}(\lambda -i0,x,y )
.$$ We have $ R_{H  _{\omega }}(\lambda -i0,x,y )=\overline{R_{H
_{\omega }}(\lambda +i0,x,y)}$,  $\overline{F_{  \pm }(x, k)}=F_{
\pm }(x, -k) $. Pick now $\lambda \in \sigma _e(H_\omega )$ with
$\lambda
>\omega $ and $x>y$. For $k\in \Bbb R\backslash S$ since $G_{ \pm
}(x, k)$ has real entries for $k\in \Bbb R$ we get
$$\aligned H(E ,x,y)&= \frac{F_{  +}(x,-k)\, \,
 {^t(}\sigma _3F_{  -}(y,-k))}{[ F_{  +}(-k) ,     F_{  -}(-k)
]}-\frac{F_{  +}(x,k)\, \, {^t(}\sigma _3F_{  -}(y, k))}{[ F_{+ }(k)
,  F_{  -}( k) ]}  \\ +&     \frac{G_{  +}(x, k) \, \, {^t(}\sigma
_3G_{  -}(y, k))}{[ G_{  +}( k) ,     G_{  -}( k) ]}-\frac{G_{
+}(x,k)\, \, {^t(}\sigma _3G_{  -}(y, k))}{[ G_{+ }(k) , G_{  -}( k)
]}
\endaligned
$$
  with the last line equal to 0. By  Lemmas 8.5-6, for $k\not \in S$

$$\aligned &  \overline{F_{
\pm } (x,k)}=F_{ \pm } (x,-k)\, , \, \overline{T(k)}=T(-k) \, , \,
\overline{R_{\pm } (k)}=R_{\pm } (-k), \\&  { F_{ + } (x,k)} =
\overline{T (k) F_{ - } (x,k) -R_{+ } (k) F_{ +} (x,k)} ,\\& F_{ - }
(y, -k)=  { {T  (k)}   {F_{ + } (y,k)} -R_{- } ( k) {F_{ -}
(y,k)}},\endaligned
$$
and so we obtain

$$\aligned  & H(E ,x,y)=   \frac{ T(k)
  \overline{F_{ + } (x,k)}\, \,  {^t(}\sigma _3  {F_{  +}(y, k)})}
  {[ F_{  +}(-k) ,     F_{  -}(-k)
]} - \frac{ \overline{T(k) }   \overline{F_{ -} (x,k)} \, \,
{^t(}\sigma _3   {F_{  -}(y, k)})}{[ F_{+  }(k) ,  F_{ -}( k) ]} \\&
-F_{ +} (x, -k) \, \,  {^t(}\sigma _3F_{  -}(y,  k)) \left [ \frac{
{R_{-}(k)}} {[ F_{+  }(-k) ,  F_{  -}(- k) ]}- \frac{
\overline{R_{+}(k)} } {[ F_{  +}( k) ,     F_{ -}( k) ] } \right ].
\endaligned $$
We claim that the last line is  zero. Set $w_0(k):=  W[ \phi
_{+}^{0}( k ), \phi _{-}^{0} ( k )]$ and multiply the bracket by
$w_0(k)$. Then we get
$$W[ \phi
_{+}^{0}(  k  ), \phi _{-}^{0} ( k )] [\dots ]= \frac{W[ \phi
_{+}^{0}( k  ), \phi _{-}^{0} ( k )] } {[ F_{+  }(-k) , F_{  -}(- k)
]}{R_{-}(k)}-\overline{R_{-}(k)} { T(k)}.$$ We have $ W[ \phi
_{+}^{0}(  k  ), \phi _{-}^{0} ( k )] =W[ \phi _{-}^{0}( - k  ),
\phi _{+}^{0} ( -k )]=-W[ \phi _{+}^{0}( - k  ), \phi _{-}^{0} ( -k
)] $, so in particular $  \overline{w} _0(k) =-w_0(k)$. Then
$$W[ \phi
_{+}^{0}(  k  ), \phi _{-}^{0} ( k )] [\dots ]= - \left ( {R_{+}(
k)} \overline{ T( k)}+\overline{R_{-}(k)} { T(k)}\right )=0.$$ This
yields the claim. We pick now $\lambda
>\omega $, $\lambda \in \sigma _e(H_\omega) $ with $\lambda =E(k)+
\omega$ and with $k\not \in S$. Picking also $x>y$ we can write  by
Lemma 5.5
$$\aligned & H(\lambda ,x,y)\frac{d\lambda }{dk} =
\dot E(k) H(\lambda ,x,y)= - i w_0(k) H(\lambda ,x,y)= \\&
   i\frac{ w_0(-k) T(k)
  \overline{F_{ + } (x,k)}\, \,  {^t(}\sigma _3  {F_{  +}(y, k)})}
  {[ F_{  +}(-k) ,    F_{  -}(-k)
]}+i\frac{ w_0(k)\overline{T(k) }   \overline{F_{ -} (x,k)} \, \,
{^t(}\sigma _3   {F_{  -}(y, k)})}{[ F_{+  }(k) ,  F_{ -}( k) ]}
\\& = i |T(k)|^2 \left (  \overline{F_{ +} (x,k)} {^t(}\sigma _3
{F_{ +}(y, k)})+ \overline{F_{ -} (x,k)}  {^t(}\sigma _3 {F_{ -}(y,
k)})\right ) \\& = 2\pi i  \left (\overline{\Psi   (x,k)}
{^t(}\sigma _3 {\Psi (y, k)})+ \overline{\Psi  (x,-k)}  {^t(}\sigma
_3 {\Psi (y, -k)})\right ).
\endaligned
$$
By continuity the formula extends to any $k> 0$. The same identity
holds for $x<y$. Hence for any $M>0$ we have
 $$\aligned & \frac{1}{2\pi i}\int  _{  \omega  }^{M^2+\omega  }
  \left [ R _{H  _{\omega }}(\lambda  -i0 )-    R_{H  _{\omega
}}(\lambda +i0 ) \right ] f  d\lambda =\\&   \int _{ 0 }^{M }dk\int
_\Bbb R dy
  \left (\overline{\Psi   (x,k)}
{^t(}\sigma _3 {\Psi (y, k)})+ \overline{\Psi  (x,-k)}  {^t(}\sigma
_3 {\Psi (y, -k)})\right ) f(y).    \endaligned
$$
Repeating the argument for $E<-\omega $ and   for $M\to \infty $ we
conclude the proof.

\proclaim{Lemma 8.13} The following operators  $P_\pm (\omega )$ are
well defined in $L^2_x$:
$$\aligned &P_+(\omega )u =\lim _{M \nearrow  \infty}\lim _{
\epsilon \searrow 0  }
 \frac 1{2\pi i}
 \int  _{[\omega ,M]\cap (\omega +\sigma
(h_0))} \left [ R_{H_\omega}(\lambda +i\epsilon )-
R_{H_\omega}(\lambda -i\epsilon )
 \right ] ud\lambda \\&
P_-(\omega )u =\lim _{M \nearrow  \infty}\lim _{\epsilon \searrow 0
}
 \frac 1{2\pi i}
 \int _{[-M,-\omega ]\cap (-\omega -\sigma
(h_0))} \left [ R_{H_\omega}(\lambda +i\epsilon )-
R_{H_\omega}(\lambda -i\epsilon ) \right ] ud\lambda
\endaligned
$$
and have kernel
$$\aligned & P_+(
 \omega )(x,y)=  \int _\Bbb R   { \overline{\Psi  (x ,k)}} \, \, \,  {^t\left
( \sigma _3  {\Psi    (y ,k)}\right ) }
   dk ,\\&
P_-(
 \omega )(x,y)= \int _\Bbb R    \sigma _1 { \overline{\Psi  (x ,k)}}\, \,  \,  {^t\left (  \sigma _3\sigma
_1  {\Psi    (y ,k)}\right ) }   dk.\endaligned $$ For any $M>0$ and
$N>0$ and for $C=C (N,M,\omega  )$ upper semicontinuous in $\omega
$, we have $$  \|  \langle x \rangle ^{M}  (P_+(\omega
 )-P_-(\omega  )-P_c(H_\omega   )\sigma _3) f\|  _{L^2_x }\le
C  \|  \langle x \rangle ^{-N}    f\|  _{L^2_x }. \tag 1$$
\endproclaim
{\it Proof.} The first two statements follow from the argument of
Lemma 8.12 and from \S 7. Now we want to prove (1).  The proof  is
similar to Lemma 5.12 \cite {C1}. For this proof we set
$H=H_\omega$, $H_0=\sigma _3 (h_0 + \omega )$,
  $R_0(z)= (H_0-z)^{-1}$ and $R(z)= (H-z)^{-1}$.
To prove (1) it is enough to write $P_c=P_++P_-$ and to prove
$\|\left [ P_\pm \sigma _3 \mp P_\pm \right ] g\|  _{L^ {2, M}_x}  \le
c\| g\| _{L^ {2, -N}_x} .$ It is not restrictive to consider only
$P_+$. Setting $H=H_0+V$, we write
$$\aligned &\sum_\pm \pm R(\lambda \pm  i\epsilon )
= \sum _\pm \pm (1+ R_0(\lambda \pm  i\epsilon ) V )^{-1}
R_0(\lambda \pm  i\epsilon ) .
\endaligned \tag 2
$$
By elementary computation
$$  R_0(\lambda \pm  i\epsilon )
\sigma _3 =R_0(\lambda \pm  i\epsilon ) -2 (h_0 +\omega +\lambda \pm
i \epsilon )^{-1} \text{diag}(0,1).$$ Therefore  $$ \text{rhs}  (2)
\sigma _3 = \text{rhs} \,(4)+2\sum _\pm \pm (1+ R_0(\lambda \pm
i\epsilon ) V )^{-1} \text{diag}(0,1)
  R_{h_0}   (- \omega -\lambda \mp i
\epsilon ) .
$$
Hence we are reduced to show that $$Ku=$$
$$\lim _{M \to  \infty}
\lim _{\epsilon \to 0^+} \sum _\pm \pm \int _{[\omega , M]\cap
\sigma (H_0)} (1+ R_0(\lambda \pm  i\epsilon ) V )^{-1}
\text{diag}(0,1)
 R_{h_0}   (- \omega -\lambda \mp i
\epsilon )
 u d\lambda
$$
defines an operator such that  for some fixed $c$  $$\| Ku \| _{L^
{2, M}_x} \le c \| u \| _{L^ {2, -N}_x} \tag 3$$
 For $m \ge 1$ we expand $(1+R_0V )^{-1}=\sum
_{j=0}^{m+1} \left [ - R_0V \right ] ^j +R_0VRV (-R_0V)^m $ and we
consider the corresponding decomposition $$K=\sum
_{j=0}^{m+1}K_j^0+\Cal K.\tag 4$$   We have $K_0^0= 0$ since for any
$u\in L^2_x$ we have
$$  \lim _{\epsilon \to 0^+}  \int
_{[\omega , M]\cap \sigma (H_0)}
 \sum _\pm \pm (h_0 +\omega +\lambda \pm i
\epsilon )^{-1} \text{diag}(0,1)ud\lambda =0.$$ We next consider
$K_j^0$  for $j\in [1,m]$ and prove $$\| K_j^0u \| _{L^ {2, M}_x} \le
c \| u \| _{L^ {2,  -N}_x}. \tag 5$$ For definiteness we consider
$j=1$, in fact the other cases are similar. We have
$$\aligned &
K_1^0u  =
 \int _{ \omega + \sigma (h_0)}  \left [  R_0 (\lambda + i
0 ) - R_0 (\lambda - i 0 ) \right ] V    R_{h_0}   (- \omega
-\lambda   )\text{diag}(0,1)
 u d\lambda .\endaligned
$$
We also have
$$K_1^0u  =\sum _{j=1}^{n_0}K_{1,j}^0u, \text{ with }
K_{1,j}^0u= \int _{\partial r _j}  R_0 (z ) V       (h_0+2 \omega +z
)^{-1}\text{diag}(0,1)
 u dz\tag 6$$
with $\{ r _j, j<n_0\}$ thin rectangles , with sides parallel to the
coordinate axes,  with each $r _j$ containing the $j$'th spectral
band  of $h_0$   in its interior and disjoint from all the other
bands. $r _{n_0}$ is  an unbounded strip containing the unbounded
spectral band in its interior. We  prove (5) for each
 $K_{1,j}^0$.
Take $j<n_0$.  Schematically
$$ K_{1,j}^0u=\int _{\partial r  _{j} } R _{h_0}(z)
V R _{h_0}   (- 2\omega -z   )  udz\tag 7$$ where now $V$ is a
scalar exponentially decreasing function and $u$ is a scalar
function. But over $\partial r  _{j}$ we have that the resolvents in
(7) have kernel in absolute value bounded by $\beta e^{-\alpha
|x-y|}$ for some fixed $\alpha >0$ and $\beta >0$. Since the length
of $\partial r _{j}$ is finite, for some $C$ and for $j<n_0$
$$|K_{1,j}^0u(x)|\le C \int _{\Bbb R^2}e^{-\alpha |x-t|}V(t)e^{-\alpha |t-y|}u(y)dtdy.
\tag 8
$$
For the latter operator one gets easily (5). Next we consider
$j=n_0$.  Once again we focus on (7). Notice that for $u$ a Schwarz
function we have that (7) for $j=n_0$ converges in   $L^{2,s}_x$ space
for $s>1/2$ because the integrand decays like $\langle \lambda
\rangle ^{-3/2}$. Using the fact that for $N\ge 1$
$$\int _{N-iN^2}^{N+iN^2} \langle z \rangle ^{-\frac{3}{2}} |dz|\le
C\langle N \rangle ^{-\frac{1}{2}},$$ we conclude that we can deform
the path integral  replacing $r_{ n_0}$ with any region $r$ defined
by $x\ge c$, $|y|\le \alpha x^2 $, with $\alpha >0$ and the
halfplane $x\ge c$ containing in its interior the unbounded band and
disjoint from the bounded bands. Then
$$K_{1,n_0}^0u=\int _{   \sigma } R _{h_0}(E(k))
V R _{h_0}   (- 2\omega -E(k)   )  u  \dot E (k) dk$$ with for $k=
 a+ib$, $\sigma$ the union of two paths of the form $b=b_0>$,  $a\in
 [\pi n(n_0), \infty )$ and $a\in
 (-\infty , -\pi n(n_0)].$ Then
$$\aligned & |K_{1,n_0}^0u(x)|\le C (A+B) \text{ with}\\&
A= \int _{|y-t|\le 1} dydt\int _{a\ge \pi n(n_0)}\frac{da}{a}
e^{-b_0|x-t|}V(t) e^{-a|y-t|}|u(y)| \\& B= \int _{|y-t|\ge 1}
dydt\int _{a\ge \pi n(n_0)}\frac{da}{a}   e^{-b_0|x-t|}V(t)
e^{-a|y-t|}|u(y)| .\endaligned$$ Then
$$\aligned & \langle x \rangle ^M A\lesssim \int _{|y-t|\le 1} dydt
\langle t \rangle ^{M+N}|V(t)|e^{-\frac{b_0}{2}|x-t|} \left | \log
(|y-t|) \right |  \langle y \rangle ^{-N}|u(y)| \\& \lesssim \| u\|
_{L_x^{2,-N}}
\endaligned$$
and   $$\aligned & \langle x \rangle ^M B\lesssim \int _{|y-t|\ge 1}
dydt \langle t \rangle ^{M+N}|V(t)|e^{-\frac{b_0}{2}|x-t|-\frac{\pi
n(n_0)}{2} |y-t|}    \langle y \rangle ^{-N}|u(y)|
\\& \lesssim \| u\| _{L_x^{2,-N}}.
\endaligned$$
This yields  (3) also for $K_{1,n_0}^0$.  Finally we focus on

$$\aligned &(-)^{m }\Cal  K u=
  \sum _\pm \\&  \pm \int _ { \omega +\sigma (h_0) }
 R_0^\pm (\lambda   ) V R^\pm (\lambda   )
V \left [ R_0^\pm (\lambda   ) V\right ] ^m \text{diag}(0,1) R_{h_0}
(-\omega -\lambda   ) ^{-1} ud\lambda \\
&=\sum _j\int _{\partial r_j }
 R_0(\zeta ) V  R(\zeta )
V\left [ R_0(\zeta ) V\right ] ^m \text{diag}(0,1)R_{h_0} (-\omega
-\zeta    )   ud\zeta .
\endaligned $$
We have $ \| V  R(\zeta ) V \| _{B(L^{2,-N}_x, L^{2, N}_x)} \lesssim
\langle z \rangle ^{-\frac{1}{2}}$ for all $N$ and by repeating the
previous argument, deforming $\partial r _{n_0} $   we obtain  also
$\| \Cal K\| _{B(L^{2,-N}_x, L^{2,M}_x)}<\infty$
  for all $M$ and $N$.

\head \S 9 Proof of Lemmas 3.1-3.2 \endhead

Lemma 3.1 is a consequence of the following statement:

\proclaim{Lemma 9.1 } Set $P_c(H_\omega )e^{itH_\omega } (x,y)=\Cal
U(t,x,y)+\Cal V(t,x,y) $ with
$$\aligned & \Cal U(t,x,y)=\int   _{|k|\le \pi n(n_0)}
e^{i t E  (k)   } \phi (x,y,k) dk \, , \\& \Cal V(t,x,y)=\int
_{|k|\ge \pi n(n_0)} e^{i t E  (k)   } \phi (x,y,k) dk\\&\text{ with
} \phi (x,y,k) ={ \overline{\Psi  (x ,k)}} \, \, \, {^t\left (
\sigma _3 {\Psi    (y ,k)}\right ) }
  +   \sigma _1 { \overline{\Psi  (x ,k)}}\,
  \,  \,  {^t\left (  \sigma _3\sigma
_1  {\Psi    (y ,k)}\right ) }  .\endaligned
 $$
Then $|\Cal U(t,x,y) | \le C         \langle t\rangle    ^{-\frac
13}$ and $|\Cal V(t,x,y) | \le C         |t|  ^{-\frac 12} .
 $
\endproclaim
{\it Proof.} Using cutoff functions we distinguish between $k$ close
to $S$ and $k$ away from $S$. In the latter case the proof is the
same of Lemma 6.5. So for $\chi (k) $ a smooth cutoff function with
small support $I$ with $I\cap S =\{  \widehat{k} \} $ and with $I
\subset \Bbb R_+$ (the above choices are not restrictive) contained
either in $|k|<\pi n(n_0)$ or in $|k|> \pi n(n_0)$, we consider the
integral
$$\int _{I}e^{i t E  (k)   }\chi (k) \overline{\Phi} _+(x,k) \, \, ^t
{\Phi} _-(y,k)\, dk .\tag 1$$  Notice that by (8.6)
$$\aligned & \Phi _{ +} (x,k)=\phi ^{0}
_\pm
 (x,k)
 \overrightarrow{e}_1  (1+A(x,k))  + \phi ^{0}
_-
 (x,k)
  B(x,k)  + C(x,y,k)
\endaligned  $$
with $|\partial _{k} ^{a}A(x,k)|+|\partial _{k}
^{a}B(x,k)|+|\partial _{k} ^{a}C(x,y,k)|< c_0$ for fixed $c_0$ for
$a=0,1.$ Then by stationary phase $|(1)|\le \langle t \rangle
^{-\alpha}$ with $\alpha  $ either $1/2$   if $I $  is  in $|k|<\pi
n(n_0)$ or  $1/3$ if $I $    in $|k|> \pi n(n_0)$.
\bigskip
Before proving Lemma 3.2 we need:

\proclaim{Lemma 9.1} Let $I_{\alpha, \beta}[f]=K_{\alpha,
\beta}(t)*f$ with $K_{\alpha, \beta}\equiv \chi_{[0,1)}t^{-\alpha} +
\chi_{[1,\infty)} t^{-\beta}$ for given $0<\alpha<1$ and $ 0<\beta
<1$. Then we have the following estimates:
$$
\|I_{\alpha, \beta}[f]\|_{\ell^p({  \Bbb {Z}}, L^r_t[n, n+1])}\leq C
\|f\|_{\ell ^q({  \Bbb {Z}}, L^s_t[n,n+1])}
$$
where $C\equiv C(p,q,r,s)>0$ and
$$\align &
1+ \frac{1}{r}\geq \alpha+\frac 1s \tag 1 \\&  1+\frac 1p\leq
\beta+\frac 1q\tag 2 \\& (r,s)\neq \left (\infty,\frac 1{1-\alpha}
\right ), \left(\frac 1\alpha, 1\right) \tag 3\\&  (p,q)\neq \left
(\infty,\frac 1{1-\beta}\right ),  \left(\frac 1\beta, 1\right).
\tag 4 \endalign
$$
If moreover we assume $\alpha=0$ and $0<\beta<1$ then
$$\|I_{0, \beta} f\|_{\ell ^p({  \Bbb {Z}}, L^r_t[n, n+1])}\leq C
\|f\|_{\ell ^q({  \Bbb {Z}}, L^s_t[n,n+1])}
$$
where $1\leq r, s\leq \infty$ and $1\leq p, q\leq \infty$ satisfy
 (2) and  (4).
\endproclaim
The elementary  proof is in Lemma 2.2 \cite{GV}.

\bigskip
{\it Proof of Lemma 3.2.} Using the definition of $\Cal U(t)$ and
$\Cal V(t)$, the proof of the Strichartz estimates for (2) and (4)
Lemma 3.2 are the same of  that in Lemma 3.1 in \cite{C1}. We turn
now to the estimate for $\Cal U(t)$. In the sequel the pairs
$(r,p)$, $(r_1,p_1)$ and $(r_2,p_2)$ are always admissible. First of
all, by $[H_\omega , \Cal U (t)]=0$, it is not restrictive to
consider only case $k=0$. We set $P_\pm =P_\pm (H _\omega )$.  We
have $[{\Cal U}  ,  P_\pm ]=0$, ${\Cal U}^\ast \sigma _3= \sigma
_3\Cal U$, $P_\pm ^\ast \sigma _3= \sigma _3P_\pm$.

{\it First step: proof of (1) Lemma 3.2.} The case $(r,p)=(\infty ,   2)$
is trivial, hence by interpolation it suffices to prove (1)
in the case $(r,p)=(4,\infty )$. By Lemma 3.1
$$\left \| \int_0^\infty {\Cal U}(t-s)   F(s)   ds
\right \|_{L^\infty_x} \leq C \int_0^\infty  \|F(s)\|_{L^1_x}
\frac{ds}{\langle t-s\rangle^\frac 13}.$$   By Lemma 9.1  with
$\alpha=0$ and $\beta=\frac 13$  $$\left \| \int_0^\infty {\Cal U}(t-s)
  F(s)    ds \right \|_{\ell ^{6}(\Bbb
Z,L^{\infty}_{t}([n,n+1],L^{\infty}_x))}
 \leq C \|F\|_{\ell^{\frac 65}
(\Bbb Z,L^{1}_{t}([n,n+1],L^{1}_x))}.\tag 9.1
$$
By Fubini theorem we have

 $$\aligned &  \left
\langle \int_0^\infty {\Cal U} (s)  P_\pm  F(s)  ds,\sigma _3 \int_0^\infty {\Cal
U}(t) P_\pm F(t) dt\right \rangle _{  x}\\& = \left \langle \int_0^\infty
{\Cal U}(t-s) P_\pm F(s)  ds, \sigma _3 P_\pm  F(t) \right\rangle _{
{t,x}}.\endaligned
$$
This implies

 $$\aligned &
\left \|\int_0^\infty {\Cal U}(s)    F(s)  ds\right \|_{L^2_x}^2\lesssim
  \left \|\int_0^\infty {\Cal U}(t-s)   F(s)
ds\right \|_{\ell ^{6}(\Bbb Z,L^{\infty}_{t}([n,n+1],L^{\infty}_x))}
\times \\&  \|   F\|_{\ell ^{\frac 65}(\Bbb
Z,L^{1}_{t}([n,n+1],L^{1}_x))}    \leq C \|F\|_{\ell ^{\frac
65}(\Bbb Z,L^{1}_{t}([n,n+1],L^{1}_x))}^2 \endaligned \tag 9.2
$$
 by (9.1) and by $\| P_\pm u \| _{L^2_x}^2 \approx \left |\langle  P_\pm u,  \sigma _3
P_\pm u  \rangle \right |$. The latter follows from the following
facts, for $W$ and $Z$ the operators in Proposition 7.1:  for $u\in
L^2_c(H_\omega )$ and  $u=Wv$, we have $\| P_+ u \| _{L^2_x} \approx
\| \text{diag} (1,0)v \| _{L^2_x}$, $\| P_- u \| _{L^2_x} \approx \|
\text{diag} (0,1)v \| _{L^2_x}$;  for any pair $\widetilde{u}=
W\widetilde{v}$ we have $\langle    \widetilde{u},  \sigma _3
\widetilde{ u}  \rangle = \langle    \widetilde{v},  \sigma _3
\widetilde{ v}  \rangle .$

   To finally    deduce (1) Lemma 3.2, notice
that by combining the Fubini theorem, a duality argument, the
H\"older inequality and (9.2),  we get for $f\in L^2_c(H_\omega )$
$$\aligned &\|{\Cal U}(t)  f \|_{\ell ^{6}(  \Bbb{Z},
L^{\infty}_{t}([n,n+1],L^{\infty}_x))} =\sup_{G\in {\Cal B}_{\frac
65, 1, 1}} \langle {\Cal U}(t)  f ,\sigma _3 G(t,x)\rangle _{ {t,x}}
\\&=\sup_{G\in {\Cal B}_{\frac 65, 1, 1}} \left
\langle  f,\sigma _3 \int_0^\infty \Cal U( t) P_c(H_\omega ) G(t)   dt\right
\rangle _{ _x}
\\&\leq \|f\|_{L^2_x}
\sup_{G\in {\Cal B}_{\frac 65, 1, 1}} \left \|  \int_0^\infty \Cal U( t)
P_c(H_\omega )
  G(t)    dt \right \|_{L^2_x} \leq C \|f\|_{L^2_x},\endaligned$$
where:
$${\Cal B}_{\frac 65, 1, 1}:=
\left \{G\in \ell^{\frac 65}(  \Bbb{Z},L^{1}_{t} ([n,n+1],L^{1}_x))
\text{ s.t. } \|G\|_{\ell ^{\frac 65}( \Bbb{Z},L^{1}_{t}
([n,n+1],L^{1}_x))} = 1 \right \}.$$

{\it Second step: proof of (3) Lemma 3.2.} We split the proof 7
subcases.

{\it First subcase:   $(r_i, p_i)=(\infty, 2)$ for $i=1,2$.}
 In this case the estimate (3) Lemma 3.2 is equivalent to
the following one:
$$\left \| \int _{0}^{t}{\Cal U}(t-s)   F(s
) ds\right \|    _{L^\infty_tL^2_x} \leq C
\|F\|_{L^1_tL^2_x}$$ whose proof is elementary:
$$\left \| \int _{0}^{t}{\Cal U}(t-s)   F(s
)  ds\right \|_{L^2_x}\leq \int _{0}^{t} \|{\Cal U}(t-s)   F(s
)\|_{L^2_x}   ds \leq  C\int _{0}^{t} \|F(s )\|_{L^2_x}   ds,$$
where we have used the fact that $\|{\Cal U}(t) f\|_{L^2_x} \leq C
\|f\|_{L^2_x}$ for a fixed $C$.

{\it Second subcase:   $(r_i, p_i)=(4, \infty)$ for $i=1,2$.}
  It is similar to the proof of  (9.1).

{\it Third subcase: $(r_1, p_1)=(4, \infty)$ and $(r_2,
p_2)=(\infty, 2)$.} In this case (3) Lemma 3.2 reduces to: $$\left
\| \int _{0}^{t} {\Cal U}(t-s)  F(s )  ds\right \|    _{\ell ^{6}(
 \Bbb{Z},L^{\infty}_{t}([n,n+1],L^{\infty}_x))}\leq C
  \|F\|_{L^1_t L^2_x}. \tag 9.3
$$
Notice that we have the identity:
$$\int _{0}^{t} {\Cal U}(t-s)   F(s
) ds= \int_0^\infty \chi(t-s){\Cal U}(t-s)  F(s )  ds,$$ with $\chi$
is the characteristic function of the half line $(0, \infty)$. By
the Minkowski inequality we get
$$\left \| \int _{0}^{t} {\Cal U}(t-s)   F(s
)  ds\right \|    _{\ell ^{6}(
\Bbb{Z},L^{\infty}_{t}([n,n+1],L^{\infty}_x))}
$$ $$\leq \int_0^\infty \|\chi(t-s) {\Cal U}(t-s)  F(s)
\|_{\ell ^{6}(  \Bbb{Z},L^{\infty}_{t}([n,n+1],L^{\infty}_x))}
  ds$$
$$\leq \int_0^\infty \|{\Cal U}(t-s)   F(s)
\|_{\ell ^{6}(  \Bbb{Z},L^{\infty}_{t}([n,n+1],L^{\infty}_x))}
\hbox{ } ds\leq C \int_0^\infty \|F(s)\|_{L^2_x}  ds$$ where we have used
  (1) Lemma 3.2 at the last step.

{\it Fourth subcase:  $(r_1, p_1)$ any admissible pair and  $(r_2,
p_2)=(\infty, 2)$.} Follows by interpolation between the first   the
third subcases.

{\it Fifth subcase:  $(r_1, p_1)=(\infty, 2)$, $(r_2, p_2)=(4,
\infty)$.} Follows by a duality argument. In fact we have
$$\left \|\int_0^t {\Cal U}(t-s)   F(s)  ds\right
\|_{L^\infty_tL^2_x}=  \sup_{G\in {\Cal B}_{1,2}} \left \langle
\int_0^t {\Cal U}(t-s)  F(s)  ds,\sigma _3 G\right\rangle _{
{t,x}}=$$
$$=\sup_{G\in {\Cal B}_{1,2}} \left \langle F,
\sigma _3 \int_t^\infty {\Cal U}(t-s)   G(s)  ds\right \rangle _{
{t,x}},
$$ where
$${\Cal B}_{1,2}:=
\left \{G\in L^1_tL^2_x \text{ s.t. } \|G\|_{L^1_tL^2_x}=1\right \}.
$$
By H\"older inequality
$$\aligned & \left \|\int_0^t {\Cal U}(t-s)   F(s)  ds\right
\|_{L^\infty_tL^2_x} \leq \|F\|_{\ell ^{\frac 65}(
\Bbb{Z},L^{1}_{t}([n,n+1],L^{1}_x))}\times \\&  \sup_{G\in {\Cal
B}_{1,2}} \left \|\int_t^\infty {\Cal U}(t-s) G(s) ds \right
\|_{\ell ^{6}(\Bbb Z,L^{\infty}_{t}([n,n+1],L^{\infty}_x))}
 \leq C \|F\|_{\ell ^{\frac 65}(  \Bbb{Z},L^{1}_{t}
([n,n+1],L^{1}_x))},\endaligned $$ where at the last step we used
the  estimate
$$\left \|\int_t^\infty {\Cal U}(t-s)   G(s) ds
\right \|_{\ell ^{6}( \Bbb{Z},L^{\infty}_{t}([n,n+1],L^{\infty}_x))}
\leq C \|G\|_{L^1_tL^2_x},
$$
whose proof is similar to (9.3).

{\it Sixth subcase:  $(r_1, p_1)$ any admissible pair and   $(r_2,
p_2)=(4, \infty)$.}
  Follows by interpolation of second and
fifth case.

{\it Seventh subcase: remaining cases.} Follows by interpolation of
fourth and   sixth cases.

\head \S 10 Extension of a result by Christ and Kiselev to
Birman-Solomjak  spaces\endhead

Given two Banach spaces   $X$ and $Y$ let $K(s,t)$  be a continuous
function with values in $B(X, Y)$. Let us introduce the operators:
$$T_K f(t)=\int_{-\infty}^\infty K(t,s)f(s) ds$$
and
$$\tilde T_K f(t) = \int_{-\infty}^t K(t,s) f(s) ds.$$
In this section we shall prove the following modified version of
Lemma 3.1 \cite{SmS}. \proclaim {Proposition 10.1} Let $1\leq p,
q,r\leq \infty$ be such that $1\leq r<\min \{p,q\}\leq \infty$.
Assume that there exist $C>0$ such that:
$$
\|T_K f\|_{\ell ^q(\Bbb Z, L^p_t([n,n+1],Y))}\leq C \|f\|_{L^r_t(X)},
\tag 10.1$$ then
$$
\|\tilde T_K f\|_{\ell ^q (\Bbb Z, L^p_t([n,n+1],Y))}\leq C'
\|f\|_{L^r_t(X)}, \tag 10.2$$ where $C'=C'(C, p,q,r)>0$ is another
suitable constant.
\endproclaim
 {\it Remark.}
In the case $p=q$ the previous proposition follows from \cite{CK}.

We shall need the following lemma.

\proclaim {Lemma 10.2} Let $k\in \N$ and $1\leq p,q,r\leq \infty$ as
in Proposition 10.1. Assume that $\{I_j\}_{j=1,...,2^k}\subset \R$
is a family of intervals (eventually unbounded) such that:
$${\R}=\cup_j I_j \text{ and }
int I_j\cap int I_k=\emptyset \hbox{ for } j\neq k.$$ Assume also
that $g_j\in \ell ^q(\Bbb Z, L^p_t([n,n+1],Y))$ for $j=1,.., 2^k$ is a
family of functions such that
$$
\|g_j\|_{\ell ^q(\Bbb Z, L^p_t([n,n+1],Y))}\leq C 2^{-\frac kr} \hbox{ }
\forall j=1,...,2^k.\tag 10.3$$ Then the following estimate holds:
$$\|\sum_{j=1}^{2^k} \chi_{I_j} g_j \|_{\ell^q(\Bbb Z, L^p_t([n,n+1],Y))}
\leq C'2^{k(\frac 1q-\frac {Min \{p, q\}}{qr})}\tag 10.4$$ where
$C'=C'(C,p,q,r)>0$ and in particular $C'$ does not depend on $g_j,
I_j, k$.
\endproclaim
  {\it Proof.} Let us fix the following notations:
$$
G(x)\equiv \sum_{j=1}^{2^k} \chi_{I_j} g_j \tag 10.5$$ and
$${\Cal I}\equiv \{I_j \hbox{ for } j=0,1,..., 2^k\},$$
where $I_j$ are the intervals given in the statement and
$I_0\equiv\emptyset$. To every $I_j\in {\Cal I}$ we associate a new
segment $\tilde I_{j}\subset I_{j}$ defined as follows:
$$\tilde I_j\equiv (z_j, z_{j+1})$$ where
$z_{j}, z_{j+1} \in {\Bbb Z}\cup \{\pm \infty\}$, $z_j\equiv \inf \{z\in
{\Bbb Z} \cup \{-\infty\}| z\in I_j\}$ and $z_{j+1}\equiv \sup \{z\in
{\Bbb Z} \cup \{+\infty\}| z\in I_j\}$. Notice that the possibility
$\tilde I_j=\emptyset$ is allowed. Next we introduce
$${\tilde {\Cal I}}\equiv \{I_j\in {\Cal I}| \tilde I_j\neq \emptyset\}
$$
and also
$$
\tilde {\Bbb Z}\equiv \{z\in {\Bbb Z}| (z, z+1)\cap (\cup_{j=1}^{2^k} \tilde
I_j)=\emptyset \}. \tag 10.6$$ Notice that in a more explicit way we
can write
$${\tilde {\Cal  I}}=\{I_{i_0}, I_{i_1},...,I_{i_h}\}$$
for suitable $0\leq h\leq 2^k$ and $0=i_0<i_1<...<i_h\leq 2^k.$
For every function $F\in l^q(\Bbb Z,
L^p_t([n,n+1],Y))$ we get:
$$\|\chi_{(\cup_{I_j\in \tilde I} \tilde I_j)} F \|_{\ell ^q(\Bbb Z, L^p_t([n,n+1],Y))}
=\left (\sum_{j=i_0,..., i_h} \| \chi_{\tilde I_j} F \|^q_{\ell
^q(\Bbb Z, L^p_t([n,n+1],Y))} \right )^\frac 1q.$$ In particular if we
choose in this identity $F=G$, with $G$ as in (10.5), then we get:
$$\aligned &
\|\chi_{(\cup_{I_j\in \tilde I} \tilde I_j)} G \|_{\ell ^q(\Bbb Z,
L^p_t([n,n+1],Y))}^q=\|\chi_{(\cup_{I_j\in \tilde I} \tilde I_j)}
(\sum_{l=1}^{2^k} \chi_{I_l} g_l )\|_{\ell ^q(\Bbb Z,
L^p_t([n,n+1],Y))}^q
\\&
= \sum_{j=i_0,i_1,..., i_h} \|\chi_{\tilde I_j} g_j \|_{\ell ^q(\Bbb Z,
L^p_t([n,n+1],Y))}^q
\\&
\leq \sum_{j=i_0,i_1,...,i_h}\|g_j\|_{\ell ^q(\Bbb Z,
L^p_t([n,n+1],Y))}^q \leq C^q h 2^{-\frac {kq}r}\endaligned \tag
10.7$$ where we have used the assumption (10.3). Next we shall
estimate
$$\|\chi_{{\R}\backslash (\cup_{I_j\in \tilde {\Cal I}}
\tilde I_j)} G \|_{l^q(\Bbb Z, L^p_t([n,n+1],Y))},$$ with $G$ as in
(10.5). In order to do that let us associate to any integer $n_0\in
\tilde {\Bbb Z}$, defined in (10.6), the following set:
$${\Cal I}_{n_0}\equiv \{I\in {\Cal I}| I\cap (n_0,n_0+1
)\neq \emptyset\}.$$
In particular we can write in a more explicit way:
$${\Cal I}_{n_0}\equiv
\{I_{i(n_0)},I_{i(n_0)+1}..., I_{i(n_0)+ h(n_0))}\}$$ where $$0\leq
h(n_0)\leq 2^k, i(n_0)\in 1,...,2^k.$$ It is easy to deduce that the
assumptions on $\{I_j\}_{j=1,..,2^k}$ imply:
$$
0\leq h\leq 2^k \text{ and } \sum_{n_0\in \tilde {\Bbb Z}}
h(n_0)\leq 2^{k+1}.\tag 10.8
$$
Next notice that for every $n_0\in \tilde {\Bbb Z}$ we have:
$$\aligned &\|\chi_{(n_0, n_0+1)} G\|_{L^p([n_0, n_0+1],Y)}^p=
\sum_{J\in {\Cal I}_{n_0}}\| \chi_{J\cap [n_0, n_0+1]}
G\|^p_{L^p([n_0, n_0+1],Y)}
\\&
=\sum_{h=0}^{h(n_0)} \|\chi_{I_{i(n_0)+h} \cap [n_0, n_0+1]}
g_{i(n_0)+h}\|_{L^p([n_0, n_0+1],Y)}^p\\& \leq \sum_{h=0}^{h(n_0)}
\|g_{i(n_0)+ h}\|_{L^p([n_0, n_0+1],Y)}^p\leq C^p (h(n_0)+1)
2^{-\frac {kp}r}\endaligned$$ and this implies
$$
\|\chi_{(n_0, n_0+1)} G\|_{L^p([n_0, n_0+1],Y)}^q \leq C^q  2^\frac
qp 2^{-\frac {kq}r} h(n_0)^\frac qp
 =C^q  2^\frac qp (2^{-\frac{kp}{r}}h(n_0))^\frac qp.\tag 10.9$$
Next we split the proof of (10.4) in two cases.

 \medskip

{\it First case: $q\geq p$}

\noindent In this case we have the following estimate:
$$
C^q  2^\frac qp (2^{-\frac{kp}{r}}h(n_0))^\frac
 qp \leq C^q 2^\frac
qp 2^{-\frac{kp}r}h(n_0), \tag 10.10$$ where at the last step we
have used the assumption $q\geq p$ and the following fact
$$2^{-\frac{kp}{r}}h(n_0)
\leq 2^{-k} h(n_0)\leq 1,$$ that in turn follows from the hypothesis
$p>r$ and $0\leq h(n_0)\leq 2^k$.

\noindent By combining (10.7) with (10.9) and (10.10) we get:
$$\aligned & \| G\|_{\ell^q(\Bbb Z, L^p_t([n,n+1],Y))}^q = \|
\chi_{(\cup_{I_j\in \tilde I} \tilde I_j)} G\|_{\ell^q(\Bbb Z,
L^p_t([n,n+1],Y))}^q + \sum_{n_0 \in \tilde {\Z}} \|G\|_{L^p([n_0,
n_0+1],Y)}^q \\& \leq C^q 2^{-\frac {kq}r} h + \sum_{n_0\in \tilde
{\Bbb Z}} C^q 2^\frac qp 2^{-\frac{kp}r}h(n_0),\endaligned
$$
that due to (10.8) implies:
$$\| G\|_{\ell ^q(\Bbb Z, L^p_t([n,n+1],Y))}^q
\leq C^q 2^{k(1-\frac {q}r)} + 2 C^q  2^\frac qp
2^{k(1-\frac{p}r)}.$$ Hence, due to the hypothesis $q\geq p>r$, we
get:
$$\| G\|_{\ell ^q(\Bbb Z, L^p_t([n,n+1],Y))}^q\leq 3 C^q  2^\frac qp
2^{k(1-\frac {p}r)}.$$

\bigskip

{\it Second case: $q<p$}

\noindent In this case we have $$ C^q  2^\frac qp
(2^{-\frac{kp}{r}}h(n_0))^\frac qp =C^q 2^\frac qp
2^{-\frac{kq}{r}}h(n_0)^\frac qp \leq C^q 2^\frac qp
2^{-\frac{kq}{r}}h(n_0)\tag 10.11$$ where at the last step we have
used $h(n_0)\geq 1$ and $q<p$. By combining (10.7) with (10.9) and
(10.11) we get:
$$\aligned & \| G\|_{\ell ^q(\Bbb Z, L^p_t([n,n+1],Y))}^q =
 \| \chi_{(\cup_{I_j\in \tilde I} \tilde I_j)}
 G\|_{\ell ^q(\Bbb Z, L^p_t([n,n+1],Y))}^q +
  \sum_{n_0 \in \tilde {\Bbb Z}} \|G\|_{L^p([n_0, n_0+1],Y)}^q \\&\leq C^q 2^{-\frac {kq}r} h + \sum_{n_0\in \tilde {\Bbb Z}}
C^q 2^\frac qp 2^{-\frac{kq}r} h(n_0)\leq 3 C^q 2^\frac qp
2^{k(1-\frac qr)},\endaligned$$ where at the last step we have used
(10.8).

 \bigskip

\noindent{\it Proof of Prop 10.1.} We shall follow the notations
used in \cite{SmS}. Let $f\in L^r_t(X)$ satisfies the following
conditions:

{\item {(1)}}$\|f\|_{L^r_t (X)}=1$; {\item {(2)}} the function
$F(t)$ defined below is a bijection between $(-\infty,\infty)$ and
$(0,1)$,
$$
F(t)=\int_{-\infty}^t \|f(s)\|_X^r ds.\tag 10.12$$ By an elementary
density argument we have that (10.2) will follow once we show that
$\|\tilde T_K f\|_{\ell ^q (\Bbb Z, L^p_t([n,n+1],Y))}\leq C'$ for any
$f\in L^r_t(X)$ that satisfies the conditions above. Next we
consider the set of all dyadic subintervals of $[0,1]$. If $I$ and
$J$ are two such subintervals, we say that $I\sim J$ if the
following conditions hold: {\item {(1)}} $I$ and $J$ have the same
lenght; {\item {(2)}} $I$ must lie on the left of $J$; {\item {(3)}}
there exist two dyadic cubes $I_0$ and $J_0$ whose lenght is twice
the lenght of $I$ and $J$, and moreover $I_0$ and $J_0$ are
adiacent.

\noindent Notice that if $J$ is fixed then there are two dyadic
intervals such that $I\sim J$. Following \cite{CK,SmS} we can write
the identity
$$\tilde T_K f=\sum_{\{I, J: I\sim J\}}
\chi_{F^{-1}J} T_K (\chi_{F^{-1}I} f),$$ where $\chi_A$ denotes the
characteristic function of the set $A$ and $F(t)$ is given in
(10.12). Due to the Minkowski inequality we get: $$\aligned &
\|\tilde T_K f\|_{\ell ^q (\Bbb Z, L^p_t([n,n+1],Y))}
 \\& \leq \sum_{k=2}^\infty
\| \sum_{\{I,J: I\sim J, |I|=2^{-k}\}} \chi_{F^{-1}J} T
(\chi_{F^{-1}I} f)\|_{\ell ^q (\Bbb Z, L^p_t([n,n+1],Y))}.
\endaligned \tag 10.13 $$
Next for every $k\geq 2$ we fix the following notations:
$$I_j^k\equiv F^{-1} \left (\frac j{2^k}, \frac{j+1}{2^k}\right )
\text{ for } j=0,1,...,2^k-1,$$
$$f_j^k= \chi_{I_{j}^k} f \hbox{ and } g_j^k=T_K(f_j^k).$$
Notice that
$$\| f_j^k\|_{L^r_t(X)}=2^{-\frac kr} \quad \forall j=0,1,..., 2^k-1$$
and due to the assumption (10.2) we deduce:
$$\aligned & \|g_j^k\|_{\ell ^q (\Bbb Z, L^p_t([n,n+1],Y))}\\&
=\| T_K (f_j^k) \|_{\ell ^q (\Bbb Z, L^p_t([n,n+1],Y))} \leq C 2^{-\frac
kr} \quad \forall j=0,..., 2^k-1.\endaligned$$ We are then in
position to use Lemma 10.2 in order to deduce:
$$\aligned & \| \sum_{\{I,J: I\sim J, |I|=2^{-k}\}}
\chi_{F^{-1}J} T_K (\chi_{F^{-1}I} f)\|_{\ell ^q (\Bbb Z,
L^p_t([n,n+1],Y))}
\\& =\| \sum_{j=2,...,2^k-1} \chi_{I_{j-1}^k} g_{j}^k(x) +
\sum_{j=2,...,2^k-1} \chi_{I_{j-2}^k} g_{j}^k(x) \|_{\ell ^q (\Bbb Z,
L^p_t([n,n+1],Y))}\\& \leq \| \sum_{j=2,...,2^k} \chi_{I_{j-1}^k}
g_{j}^k(x)\|_{\ell ^q (\Bbb Z, L^p_t([n,n+1],Y))} +\|\sum_{j=2,...,2^k}
\chi_{I_{j-2}^k} g_{j}^k(x) \|_{\ell ^q (\Bbb Z, L^p_t([n,n+1],Y))}\\&
\leq 2 C' 2^{k(\frac 1q-\frac {Min\{p,q\}}{qr})}.\endaligned$$ By
(10.13) we finally get
$$\|\tilde T_K f\|_{\ell ^q (\Bbb Z, L^p_t([n,n+1],Y))}
\leq 2C' \sum_{k=2}^\infty  2^{k(\frac 1q-\frac {Min\{p,q\}}{qr})}$$
and since $\frac 1q-\frac {Min\{p,q\}}{qr}<0$ we have the
desired result.

\enddocument